\documentclass[
submission
]{dmtcs-episciences}

\usepackage[utf8]{inputenc}
\usepackage{subfigure}

\usepackage[round]{natbib}

\usepackage{amsmath}

\newtheorem{theorem}{Theorem}[section]
\newtheorem{lemma}[theorem]{Lemma}

\newtheorem{problem}[theorem]{Problem}
\newtheorem{observation}[theorem]{Observation}

\newtheorem{remark}[theorem]{Remark}

\newcommand{\nhalb}{\ensuremath{\left\lfloor\frac{n}{2}\right\rfloor}}
\newcommand{\seqnum}[1]{{\rm \underline{#1}}}

\author{Stephan Dominique Andres\affiliationmark{1}\and Igor Re\v{z}an\affiliationmark{2}}
\title{Counting Graphs Admitting a Strictly Matched Involution}
\affiliation{
  Institute of Mathematics and Computer Science, University of Greifswald, Germany\\
  Department of Mathematics and Computer Science, FernUniversit\"{a}t in Hagen, Germany} 
\accepted{}
\begin{document}
\publicationdata
{Preprint submitted to DMTCS}{2026}{1}{10.46298/dmtcs.10472}
{2026-11-04}
{None}
\maketitle
\begin{abstract}
We study \emph{graphs admitting a strictly matched involution},
\textit{i.e.}, graphs for which a strictly matched involution exists.
A \emph{strictly matched involution} of a graph $G$ 
is a self-inverse graph automorphism $f$
with the properties that the fixed point set of $f$ induces a complete subgraph in $G$
and, for the vertices that are no fixed points, every vertex is connected to
its image under $f$ by an edge. 
Andres, Huggan, Mc Inerney, and Nowakowski (2019) 
determined
upper and
lower bounds for the number $A(n)$ of 
graphs admitting a strictly matched involution
on $n$ vertices. 
We improve the upper bound and we compute the 
exact values
for the number $A(n)$ for all \begin{math}n\leq11\end{math}.
Furthermore, 
we use combinatorial enumeration techniques and a decomposition lemma to
determine the number $C(n,m)$ of 
graphs admitting a strictly matched involution
with $n$ vertices
and $m$ edges for several infinite classes of pairs $(n,m)$.
\end{abstract}

\section{Introduction}

Let \begin{math}G=(V,E)\end{math} be a (simple finite) graph. 
A \emph{graph automorphism} of $G$ is a bijective map \begin{math}\sigma:V\longrightarrow V\end{math}
with the property that, for every two vertices \begin{math}v,w\in V\end{math}, the vertices $v$ and $w$
are connected by an edge if, and only if, the vertices \begin{math}\sigma(v)\end{math} and \begin{math}\sigma(w)\end{math}
are connected by an edge. 
An \emph{involution} of $G$ is a graph automorphism $\sigma$
with \begin{math}\sigma(\sigma(v))=v\end{math} for every vertex \begin{math}v\in V\end{math}. 
By an involution~$\sigma$ of~$G$,
the vertex set of $G$ is separated into two subsets: the set 
\begin{displaymath}F_1(G,\sigma):=\{v\in V\mid \sigma(v)=v\}\end{displaymath}
of \emph{fixed points} and the set
\begin{displaymath}F_2(G,\sigma):=V\setminus F_1(G,\sigma)\end{displaymath}
of points of orbit~2. An involution $\sigma$ of~$G$ is called
\emph{strictly matched} if
\begin{enumerate}[(1)]
\item \begin{math}\forall v\in F_2(G,\sigma): v\sigma(v)\in E\end{math} and
\item \begin{math}F_1(G,\sigma)\end{math} induces a \emph{clique}, 
\textit{i.e.}, a graph in which every two vertices are connected by an edge.
\end{enumerate}
Then, an edge of type \begin{math}v\sigma(v)\end{math} is called a \emph{matching edge} of the pair \begin{math}(G,\sigma)\end{math}.

A \emph{graph admitting a strictly matched involution}
is a graph for which a strictly matched involution exists.

\pagebreak[3]

\cite{OCG} 
introduced graphs admitting a strictly matched involution, 
motivated by a specific game, the
\emph{orthogonal colouring game}. In this context, 
they
proved that,
if graphs
admitting a strictly matched involution are the boards of 
the orthogonal colouring game,
then the second player always has a strategy to force a draw \cite[Thm 1]{OCG}. 
The orthogonal colouring game is a representative of a large class of games,
the so-called \emph{scoring games},
which recently have received considerable attention in game theory by 
several authors, \textit{e.g.}, by
\cite{SS}, \cite{johnson}, \cite{La}, \cite{Lb}, \cite{Lc}, 
\cite{Ld}, 
\cite{stewarts}, 
\cite{RLP}, 
\cite{davies}, or \cite{Lf}.
Graphs admitting a strictly matched involution are also interesting
from a computational perspective:
\cite{complexityOCG} established the NP-completeness of the decision problem
to recognise whether a graph admits a strictly matched involution. 

In the present paper, we focus on counting (the isomorphism classes of) graphs admitting a
strictly involution with $n$ vertices and $m$ edges for 
\begin{math}m,n\in{\mathbb N}\end{math}.
Moreover, for \begin{math}m,n\in{\mathbb N}\end{math}, let
\begin{itemize}
\item
\begin{math}A(n)\end{math} denote the number of graphs admitting a strictly matched involution
with $n$ vertices,
\item
\begin{math}B(m)\end{math} denote the number of graphs admitting a strictly matched involution
with $m$ edges, and
\item
\begin{math}C(n,m)\end{math} denote the number of graphs admitting a strictly matched involution
with $n$ vertices and $m$ edges,
\end{itemize}
respectively, up to isomorphism. We aim to obtain partial answers to the following problems.

\begin{problem}\label{prob:countA}
Given \begin{math}n\in{\mathbb N}\end{math} as input, determine $A(n)$.
\end{problem}

\begin{problem}\label{prob:countB}
Given \begin{math}m\in{\mathbb N}\end{math} as input, determine $B(m)$.
\end{problem}

\begin{problem}\label{prob:countC}
Given \begin{math}n,m\in{\mathbb N}\end{math} as input, determine $C(n,m)$.
\end{problem}

As main results we determine the exact values $C(n,m)$ for several infinite
classes of pairs $(n,m)$.
Furthermore, we improve a general upper bound given by \cite{OCG}
for the numbers $A(n)$.
Finally, we implement a computer program that determines $A(n)$ and $C(n,m)$
for $0\le n \le 11$.

The structure of the paper is as follows.
After fixing some notation in Section~\ref{sec:notation}
and a discussion of previous results related to this paper
in Sections~\ref{sec:XXA} and~\ref{sec:XXB},
in Section~\ref{sec:eigen} we list our main results.
Section~\ref{sec:proof}
contains the proof of the new upper bound (Theorem~\ref{oberneu}). 
In Section~\ref{sec:ex} we present the results of the computer program. 
Sections~\ref{sec:exAB} and~\ref{sec:exC} are devoted to the proofs of 
our mathematical main results (Observation~\ref{obs:A},
Theorem~\ref{obs:B},  and 
Theorem~\ref{obs:C}). We conclude with some open problems in 
Section~\ref{sec:open}.
The idea and the source code of the program are displayed in Appendix \ref{sec:append}.

\section{Preliminaries}

\subsection{Notation}
\label{sec:notation}

Let \begin{math}{\mathbb N}\end{math} denote the set \begin{math}\{0,1,2,3,\ldots\}\end{math} of nonnegative integers.
For \begin{math}k\in{\mathbb N}\end{math}, let \begin{math}\binom{V}{k}\end{math} denote the
set of all subsets of \begin{math}V\end{math} with exactly \begin{math}k\end{math} elements.
A \emph{graph} is a pair \begin{math}(V,E)\end{math} consisting of a finite set $V$,
the set of \emph{vertices}, and a set \begin{math}E\subseteq\binom{V}{2}\end{math},
the set of \emph{edges}.
By this definition, all graphs in this paper are assumed to be simple and finite.
Edges \begin{math}\{v,w\}\in E\end{math} are usually written as~\begin{math}vw\end{math}.
Two graphs \begin{math}G=(V,E)\end{math} and \begin{math}G'=(V',E')\end{math} are \emph{isomorphic}
if there is a bijective map \begin{math}f:V\longrightarrow V'\end{math} satisfying
\begin{displaymath}\forall v,w\in V:\quad\left(vw\in E\Longleftrightarrow f(v)f(w)\in E'\right).\end{displaymath}\begin{figure}[htbp]
\begin{center}
\begin{minipage}[b]{0.22\textwidth}
\begin{center}
\includegraphics[scale=0.45]{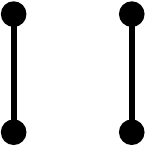}

\begin{math}2K_2\end{math}
\end{center}
\end{minipage}
\begin{minipage}[b]{0.22\textwidth}
\begin{center}
\includegraphics[scale=0.45]{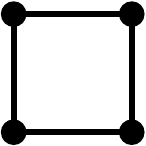}

\begin{math}C_4\end{math}
\end{center}
\end{minipage}
\begin{minipage}[b]{0.22\textwidth}
\begin{center}
\includegraphics[scale=0.45]{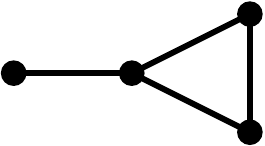}

\begin{math}P\end{math}
\end{center}
\end{minipage}
\begin{minipage}[b]{0.22\textwidth}
\begin{center}
\includegraphics[scale=0.45]{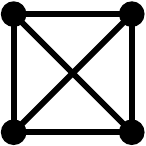}

\begin{math}K_4\end{math}
\end{center}
\end{minipage}
\end{center}
\caption{\label{fig:fourFIG}The graphs admitting a strictly matched involution
with 4 vertices}
\end{figure}Such a map is called a \emph{graph isomorphism} from $G$ to~$G'$.
Following standard graph theory conventions,
isomorphic graphs are considered as the same objects.
A \emph{graph automorphism} of $G$ is a graph isomorphism from $G$ to~$G$.

For \begin{math}n\in{\mathbb N}\end{math}, let \begin{math}K_n\end{math} denote 
the \emph{clique} (or also called \emph{complete graph}) on $n$ 
vertices, 
\textit{i.e.}, 
\begin{displaymath}K_n=\left(V,\binom{V}{2}\right)\end{displaymath} 
with \begin{math}V=\{1,\ldots,n\}\end{math}. 
Some special graphs admitting a strictly matched involution
with 4 vertices are the 4-cycle, denoted by \begin{math}C_4\end{math}, and the \emph{paw} 
(sometimes also called \emph{3-pan}), which is denoted by $P$ in this paper,
see Figure~\ref{fig:fourFIG}. The \emph{disjoint union} of two graphs $G_1$ and $G_2$,
\textit{i.e.}, a graph \begin{math}(V_1'\cup V_2',E_1'\cup E_2')\end{math} for an isomorphic copy
\begin{math}(V_1',E_1')\end{math} of $G_1$ and an isomorphic copy 
\begin{math}(V_2',E_2')\end{math} of $G_2$ with \begin{math}V_1'\cap V_2'=\emptyset\end{math},
is denoted by \begin{math}G_1\cup G_2\end{math}; and, for \begin{math}n\in{\mathbb N}\end{math}, 
the \emph{$n$-fold disjoint union} \begin{math}nG_1\end{math} is defined recursively
by \begin{math}0G_1:=(\emptyset,\emptyset)\end{math} and 
\begin{math}nG_1:=((n-1)G_1)\cup G_1\end{math} for \begin{math}n\in{\mathbb N}\setminus\{0\}\end{math}.
The \emph{complete join} of two graphs $G_1$ and $G_2$, \textit{i.e.}, the graph that 
is generated from the 
disjoint union \begin{math}G_1\cup G_2\end{math} by connecting every vertex of the copy of $G_1$
with every vertex of the copy of $G_2$, is denoted by \begin{math}G_1\vee G_2\end{math}.
For a graph \begin{math}G=(V,E)\end{math}, let \begin{math}\overline{G}\end{math} denote the \emph{complement}
\begin{displaymath}\left(V,\binom{V}{2}\setminus E\right)\end{displaymath} 
of~$G$, and, if \begin{math}E'\subseteq E\end{math}, then
\begin{math}G-E'\end{math} is the graph \begin{math}(V,E\setminus E')\end{math}. 
If \begin{math}v\in V\end{math}, then \begin{math}G-v\end{math} is the subgraph of $G$ induced by \begin{math}V\setminus\{v\}\end{math},
\textit{i.e.}, the graph obtained from $G$ by deleting the vertex $v$ and all edges
containing~$v$. The \emph{degree} \begin{math}d_G(v)\end{math} of a vertex $v$ 
in 
$G$ is the number of the
edges containing~$v$.

\subsection{The structure of graphs admitting a strictly matched involution}
\label{sec:XXA}

Our analysis is based on the following result of \cite{OCG}, which
describes the structure of graphs admitting a strictly matched involution.

\begin{theorem}[\cite{OCG}, Cor~9]
\label{thm:structure}
Every graph \begin{math}G=(V(G),E(G))\end{math} on $n$ vertices admitting a strictly matched involution
has a partition of its vertex set into three (possibly empty) vertex subsets 
inducing a clique \begin{math}C=(V(C),E(C))\end{math} 
of size \begin{math}n-2k\end{math} and two isomorphic graphs 
\begin{math}H=(V(H),E(H))\end{math} and \begin{math}H'=(V(H'),E(H'))\end{math}, 
each of size~$k$, for some \begin{math}k\in{\mathbb N}\end{math} with \begin{math}0\leq k\leq\left\lfloor\frac{n}{2}\right\rfloor\end{math}, respectively.
Moreover, there is a graph isomorphism \begin{math}\tilde{\sigma}:V(H)\longrightarrow V(H')\end{math} such that
\begin{enumerate}[(1)]
\item
\label{thm:stra}for every pair \begin{math}(v,v')\end{math} of vertices \begin{math}v\in V(H)\end{math} and \begin{math}v'\in V(H')\end{math}
with \begin{math}v'=\tilde{\sigma}(v)\end{math}
and every vertex \begin{math}x\in V(C)\end{math},
either both or none of \begin{math}vx\end{math} and \begin{math}v'x\end{math} 
are edges in \begin{math}E(G)\end{math};
\item
\label{thm:strb}for every pair \begin{math}(v,v')\end{math} of vertices \begin{math}v\in V(H)\end{math} and \begin{math}v'\in V(H')\end{math}
with \begin{math}v'=\tilde{\sigma}(v)\end{math}, we have 
\begin{math}vv'\in E(G)\end{math};
\item
\label{thm:strc}for every two pairs \begin{math}(v,v')\end{math} and \begin{math}(w,w')\end{math} of vertices 
with \begin{math}v,w\in V(H)\end{math} and \begin{math}v',w'\in V(H')\end{math},
\begin{math}v'=\tilde{\sigma}(v)\end{math}, and \begin{math}w'=\tilde{\sigma}(w)\end{math}, 
either both or none of \begin{math}vw'\end{math} and \begin{math}v'w\end{math} 
are edges in \begin{math}E(G)\end{math}.
\end{enumerate}
\end{theorem}

\pagebreak[3]

\begin{figure}[htbp]
\begin{center}
\includegraphics[scale=0.44]{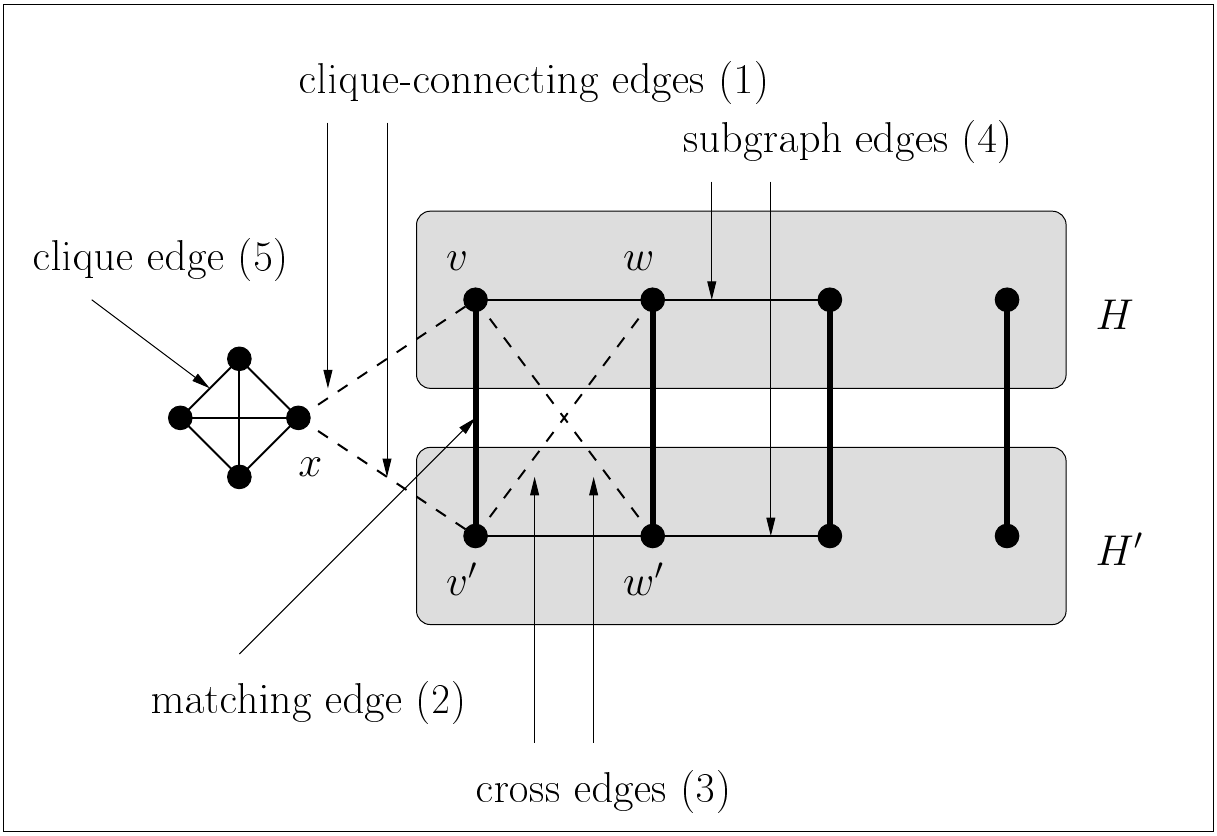}
\end{center}
\caption{\label{fig:structure}Structure of graphs
admitting a strictly matched involution~$\sigma$
according to Theorem~\ref{thm:structure} and the
five different types of edges of such a graph
with regard to~$\sigma$. The edges of type
(1), (2), and~(3) refer to (\ref{thm:stra}), (\ref{thm:strb}), and~(\ref{thm:strc}) in
Theorem~\ref{thm:structure}.}
\end{figure}

We remark that the graph isomorphism \begin{math}\tilde{\sigma}\end{math}
between the subgraphs $H$ and $H'$ 
in Theorem~\ref{thm:structure} can be extended to a 
strictly matched involution \begin{math}\sigma\end{math} of the whole
graph $G$ by defining
\begin{displaymath}
\sigma(v)=\left\{\begin{array}{ll}
\tilde{\sigma}(v)&\text{if}\ v\in V(H);\\
\tilde{\sigma}^{-1}(v)&\text{if}\ v\in V(H');\\
v&\text{if}\ v\in C.
\end{array}\right.
\end{displaymath}

As illustrated in Figure~\ref{fig:structure},
by Theorem~\ref{thm:structure}, the edge set of a graph admitting a strictly matched involution
can be decomposed into
\begin{enumerate}[(1)]
\item pairs of \emph{clique-connecting edges}, \textit{i.e.}, edges
implied by Theorem~\ref{thm:structure}~(\ref{thm:stra});
\item matching edges (cf.\ Theorem~\ref{thm:structure}~(\ref{thm:strb}));
\item pairs of \emph{cross edges}, \textit{i.e.}, the edges implied by
Theorem~\ref{thm:structure}~(\ref{thm:strc});
\item pairs of \emph{subgraph edges}, \textit{i.e.}, the edges of $H$ and $H'$; and
\item \emph{clique edges}, \textit{i.e.}, the edges induced by $C$.
\end{enumerate}

However, there might be many ways of such a decomposition, depending on the choice of $H$, $C$, and the matching edges.
For example, in the worst case, \textit{i.e.}, if $G$ is the complete graph $K_n$ on $n$ vertices, there are
\begin{displaymath}d(n)=\sum_{k=0}^{\nhalb}\binom{n}{2k}2^k\prod_{j=1}^k(2j-1)
=\sum_{k=0}^{\nhalb}\binom{n}{2k}\binom{2k}{k}k!
=\sum_{k=0}^{\nhalb}\frac{n!}{k!(n-2k)!}
\end{displaymath}
ways to obtain such a decomposition \begin{math}(C,H,H')\end{math} of the vertex set of $G$ together with
a graph isomorphism \begin{math}\tilde{\sigma}:H\longrightarrow H'\end{math}
as 
\begin{itemize}
\item
\begin{math}\binom{n}{2k}=\binom{n}{n-2k}\end{math}
is the number of ways to choose \begin{math}n-2k\end{math} vertices for the clique $C$,
\item 
\begin{math}\prod_{j=1}^k(2j-1)\end{math}
is the number 
of perfect matchings in a \begin{math}K_{2k}\end{math}, and 
\item
there are \begin{math}2^k\end{math} ways in total
to choose which
vertex of a matching edge belongs to $H$. 
\end{itemize}
We remark that the numbers $d(n)$ constitute the sequence
\seqnum{A047974} in The On-Line Encyclopedia of Integer Sequences (\cite{oeis}).

The numbers of this sequence occur in many different 
mathematical domains and applications, \textit{e.g.},
in a certain model of ad hoc mobile networks
described by \cite{LT},
in the theory of operator algebras used in
mathematical physics
(cf.\ \cite{amdeberhanetal}),
as numbers of palindromic assembly words 
(\textit{i.e.}, double occurrence words that are in a
certain way isomorphic to its inverse word), 
which were investigated by
\cite{burnsetal}
in the context of DNA recombination,
as total numbers of symmetric diagrammes in a Breuer
algebra in the theory of diagram algebras
(cf.\ \cite{HR}),
or as numbers of posets on $n$ elements
the components of which consist of at most 2 elements 
(cf.\ \cite{aljohanietal}).
Furthermore, \cite{bala} remarked 
a nice property of the
sequence \seqnum{A047974}: 
for any $k\in\naturals\setminus\{0\}$, 
it is $k$-periodic modulo $k$.
Some more general and more complicated
counting sequences that are related to the
sequence \seqnum{A047974}
are given, \textit{e.g.}, by \cite{QK}.

We consider 
counting graphs admitting a strictly matched involution, which is
another, more complicated way of counting
related to the sequence \seqnum{A047974}: starting
with the decompositions $(C,H,H')$ of the complete
graph $K_n$ with respect to a graph isomorphism 
\begin{math}\tilde{\sigma}\end{math}
from $H$ to $H'$ (which are the objects counted by the
sequence \seqnum{A047974}), for each way of deleting pairs
of clique-connecting edges, subgraph edges, and/or
cross edges, 
we obtain a graph admitting a strictly matched involution
and, 
by Theorem~\ref{thm:structure}, 
we obtain every graph admitting a strictly matched
involution on $n$ vertices by this procedure. However, at this moment,
we have constructed many isomorphic graphs. In order
to count them, we have to group them into isomorphism
classes. 
As graphs admitting a strictly matched involution
can be derived in this very natural way from the
combinatorial structures 
defining the sequence \seqnum{A047974}, which is of
significant interest in many mathematical domains, the study
of graphs admitting a strictly matched involution
(besides its importance in game theory)
might turn out to be of general interest, too. 

\subsection{Counting graphs admitting a strictly matched involution}
\label{sec:XXB}

For the number $A(n)$ of (isomorphism classes of) 
graphs admitting a strictly matched involution with
$n$ vertices, \textit{i.e.}, the entries of the sequence \seqnum{A363771} in the OEIS (\cite{oeis}), 
the following lower and upper bounds (with
a large gap inbetween) are known. The bounds depend on the number $g(n)$
of all (isomorphism classes of) graphs with $n$ vertices, 
\textit{i.e.}, the entries of the sequence \seqnum{A000088} in the OEIS (\cite{oeis}).

\begin{theorem}[\cite{errOCG}, Thm 2]\label{unter}
For every \begin{math}n\in{\mathbb{N}}\end{math},
\begin{displaymath}A(n)\geq g\left(\left\lfloor\frac{n}{2}\right\rfloor\right).\end{displaymath}
\end{theorem}

\begin{theorem}[\cite{OCG}, Thm 11]\label{ober}
For every \begin{math}n\in{\mathbb{N}}\end{math},
\begin{displaymath}A(n)\leq\sum_{k=0}^{\left\lfloor\frac{n}{2}\right\rfloor}g(k)2^{\binom{k}{2}}2^{(n-2k)k}.\end{displaymath}
\end{theorem}

\section{Our counting results}
\label{sec:eigen}

The first result of the present paper is the following improvement of
Theorem~\ref{ober}. 

\begin{theorem}\label{oberneu}
For every \begin{math}n\in{\mathbb{N}}\end{math},
\begin{displaymath}A(n)\leq\sum_{k=0}^{\left\lfloor\frac{n}{2}\right\rfloor}g(k)
2^{\binom{k}{2}}\binom{2^k+n-2k-1}{n-2k}.\end{displaymath}
\end{theorem}

Moreover, we use the construction of Theorem~\ref{oberneu} to
develop a computer program that enumerates and counts the 
number \begin{math}C(n,m)\end{math} of (isomorphism classes of) graphs
admitting a strictly matched involution with $n$ vertices and $m$ edges.
Thereby, we obtain all values of \begin{math}C(n,m)\end{math} with \begin{math}0\leq n\leq11\end{math}.

Furthermore, we determine the following exact values for \begin{math}C(n,m)\end{math}
for some pairs \begin{math}(n,m)\end{math}, as given in 
Observation~\ref{obs:A}, Theorem~\ref{obs:B}, and Theorem~\ref{obs:C}.

\begin{observation}\label{obs:A}
Let \begin{math}m,n\in{\mathbb N}\end{math}. Then
\begin{enumerate}[(a)]
\item \label{oAa}\begin{math}C(n,m)=0\end{math} if \begin{math}m<\nhalb\end{math};
\item \label{oAb}\begin{math}C(n,m)=0\end{math} if \begin{math}m>\binom{n}{2}\end{math};
\item \label{oAc}\begin{math}C(n,m)=0\end{math} if \begin{math}m\equiv0\pmod{2}\end{math} and \begin{math}n\equiv2,3\pmod{4}\end{math};
\item \label{oAd}\begin{math}C(n,m)=0\end{math} if \begin{math}m\equiv1\pmod{2}\end{math} and \begin{math}n\equiv0,1\pmod{4}\end{math}.
\end{enumerate}
\end{observation}

\begin{theorem}\label{obs:B}
Let \begin{math}n\in{\mathbb N}\end{math}. Then
\begin{enumerate}[(a)]
\item \label{oBa}\begin{math}C(n,m)=1\end{math} if \begin{math}m=\nhalb\end{math};
\item \label{oBb}\begin{math}C(n,m)=2\end{math} if \begin{math}m=\nhalb+2\end{math} and \begin{math}n\geq4\end{math};
\item \label{oBc}\begin{math}C(n,m)=7\end{math} if \begin{math}m=\nhalb+4\end{math}, \begin{math}n\equiv0\pmod{2}\end{math}, and \begin{math}n\geq8\end{math};
\item \label{oBd}\begin{math}C(n,m)=6\end{math} if \begin{math}m=\nhalb+4\end{math}, \begin{math}n\equiv1\pmod{2}\end{math}, and \begin{math}n\geq8\end{math}.
\end{enumerate}
\end{theorem}

\begin{theorem}\label{obs:C}
Let \begin{math}n\in{\mathbb N}\end{math}. Then
\begin{enumerate}[(a)]
\item \label{oCa}\begin{math}C(n,m)=1\end{math} if \begin{math}m=\binom{n}{2}\end{math};
\item \label{oCb}\begin{math}C(n,m)=2\end{math} if \begin{math}m=\binom{n}{2}-2\end{math} and \begin{math}n\geq4\end{math};
\item \label{oCc}\begin{math}C(n,m)=6\end{math} if \begin{math}m=\binom{n}{2}-4\end{math} and \begin{math}n\geq8\end{math}.
\end{enumerate}
\end{theorem}

Using the results of Observation~\ref{obs:A} and Theorem~\ref{obs:B} and the
results of the computer program, we are able to determine the number
$B(m)$ of (isomorphism classes of) graphs admitting a strictly matched
involution with $m$ edges, \textit{i.e.}, the entries of the sequence \seqnum{A363772} in the OEIS (\cite{oeis}))
for \begin{math}0\leq m\leq11\end{math}.

\section{Proof of the improved upper bound}
\label{sec:proof}

\begin{proof}[of Theorem~\ref{oberneu}]
Let $G$ be a graph admitting a strictly matched involution with $n$ vertices.
By Theorem~\ref{thm:structure}, the vertex set of $G$ can be partitioned in
two subsets $V$ and $V'$, each of cardinality~$k$ for some \begin{math}k\in{\mathbb{N}}\end{math}
with \begin{math}0\leq k\leq\left\lfloor\frac{n}{2}\right\rfloor\end{math}, and a subset $W$ of
cardinality \begin{math}n-2k\end{math} such that $W$ induces a complete graph \begin{math}C=\left(W,\binom{W}{2}\right)\end{math}
and $V$ and $V'$ induce two isomorphic graphs \begin{math}H=(V,E)\end{math} and \begin{math}H'=(V',E')\end{math}.
By property (2) in Theorem~\ref{thm:structure}, 
vertices that are mapped to each other by the isomorphism are
connected by matching edges. We have the following choices.

\begin{enumerate}[(1)]
\item \label{oberneuA}
We have \begin{math}g(k)\end{math} many possibilities to choose the edges in~$E$.
\item \label{oberneuB}
For each of the \begin{math}\binom{k}{2}\end{math} unordered pairs \begin{math}\{(v,v'),(w,w')\}\end{math} 
of pairs \begin{math}(v,v')\end{math} and \begin{math}(w,w')\end{math} of corresponding matched vertices \begin{math}v,w\in V\end{math} and
\begin{math}v',w'\in V'\end{math}, by property (\ref{thm:strc}) in Theorem~\ref{thm:structure},
we have two choices: either the crossing edges \begin{math}vw'\end{math} and \begin{math}v'w\end{math} are both in $G$
or none of them; this gives us 
\begin{math}2^{\binom{k}{2}}\end{math}
choices. 
\item \label{oberneuC}
By property (\ref{thm:stra}) in Theorem~\ref{thm:structure}, 
each of the \begin{math}n-2k\end{math} vertices in the clique $C$ is adjacent to 
each vertex of
one of the \begin{math}2^k\end{math} subsets
of the partition $V$. Since the vertices of the clique are indistinguishable, this is equivalent to
putting \begin{math}K:=n-2k\end{math} indistinguishable balls in \begin{math}N:=2^k\end{math} numbered cases (possibly several balls in a case).
Thus we get as many choices as the number of $K$-combinations of $N$ with repetition, which is
\begin{displaymath}\binom{N+K-1}{K}=\binom{2^k+n-2k-1}{n-2k}.\end{displaymath}
\end{enumerate}
As the three choices are independent, we get in total
\begin{displaymath}g(k)2^{\binom{k}{2}}\binom{2^k+n-2k-1}{n-2k}\end{displaymath}
ways to choose the edges for fixed~$k$,
and summing over all possible $k$ gives the upper bound for \begin{math}A(n)\end{math} in the theorem. 
\end{proof}

In the proof of Theorem~\ref{ober}, 
\cite{OCG}
basically used the same idea as in 
the above proof of Theorem~\ref{oberneu},
but they 
assumed that the vertices of the cliques are distinguishable, 
which gives us the number 
of $K$-variations of~$N$
(for \begin{math}K=n-2k\end{math} and \begin{math}N=2^k\end{math}) 
\begin{displaymath}N^K=(2^k)^{n-2k}=2^{(n-2k)k}\end{displaymath}
instead of the number of $K$-combinations of $N$. As the number of $K$-variations of $N$ is an upper bound for
the number of $K$-combinations of $N$ we have
\begin{displaymath}\binom{2^k+n-2k-1}{n-2k}\leq2^{(n-2k)k},\end{displaymath}
thus the new upper bound from Theorem~\ref{oberneu} is better than the old upper bound from Theorem~\ref{ober}.

\section{Exact values: computational results} 
\label{sec:ex}

\cite{rezan} computed the values of \begin{math}A(n)\end{math} for \begin{math}0\leq n\leq10\end{math}. 
Theorem~\ref{oberneu} gives an improved idea how to enumerate the graphs
admitting a strictly matched involution on $n$ vertices by computer:
we construct all configurations that are counted in the upper bound
of Theorem~\ref{oberneu} and then make isomorphism tests to exclude
isomorphic copies. The program 
based on the improved idea
is displayed in Appendix~\ref{sec:appendB}.

By this approach, we compute the exact values of the sequence \begin{math}(A(n))_{n\in{\mathbb{N}}}\end{math} 
for all \begin{math}n\in{\mathbb{N}}\end{math} in the range \begin{math}0\leq n\leq11\end{math}, thus, being able to deal with the case \begin{math}n=11\end{math}, too.
The running time for \begin{math}0\leq n\leq9\end{math} on a standard laptop 
was in 38 seconds, the case \begin{math}n=10\end{math} needed
27 minutes, and the case \begin{math}n=11\end{math} more than 31 hours.

\begin{table}[htbp]
\begin{displaymath}\begin{array}{|r|r|r|r|r|}
\hline
&&&&\\
n&\mbox{\small\begin{math}g\left(\left\lfloor\frac{n}{2}\right\rfloor\right)\end{math}}
&A(n)
&\mbox{\small\begin{math}\sum\limits_{k=0}^{\left\lfloor\frac{n}{2}\right\rfloor}g(k)2^{\binom{k}{2}}\binom{2^k+n-2k-1}{n-2k}\end{math}}
&\mbox{\small\begin{math}\sum\limits_{k=0}^{\left\lfloor\frac{n}{2}\right\rfloor}g(k)2^{\binom{k}{2}}2^{(n-2k)k}\end{math}}\\
&&&&\\
\hline
&&&&\\
0&1&1&1&1\\
1&1&1&1&1\\
2&1&1&2&2\\
3&1&2&3&3\\
4&2&4&8&9\\
5&2&9&21&25\\
6&4&21&78&113\\
7&4&65&343&545\\
8&11&240&2004&3841\\
9&11&1128&15337&31873\\
10&34&6764&141466&362753\\
11&34&53971&1714411&5112321\\
\hline
\end{array}\end{displaymath}
\caption{\label{tab:an}Lower bound (second column), exact value (third column), new upper bound (fourth column), 
and previously known upper bound (fifth column), for
the sequence $A(n)$}
\end{table}

We implemented the values $A(n)$ in the OEIS (\cite{oeis}) in sequence \seqnum{A363771}. 
In Table~\ref{tab:an}
we list those values together with the values for the lower bound from Theorem~\ref{unter} and the
upper bounds from Theorem~\ref{oberneu} and Theorem~\ref{ober}. We remark that the slight improvement
given in Theorem~\ref{oberneu} is significant for large~$n$. However, there are still large
significant gaps between the lower and new upper bound and the exact values of $A(n)$ for large~$n$.
Thus further improvements of the bounds might be an interesting topic of further research.

Moreover, with the program we indeed determined the values for \begin{math}C(n,m)\end{math}
for all \begin{math}0\leq n\leq11\end{math} and \begin{math}0\leq m\leq\binom{n}{2}\end{math}. 
The results are displayed in Table~\ref{Anm}.

\begin{figure}[htbp]
\begin{center}
\includegraphics[scale=0.25]{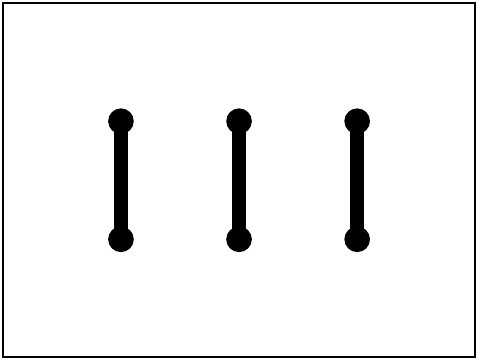}

\vspace{0.5ex}

\includegraphics[scale=0.25]{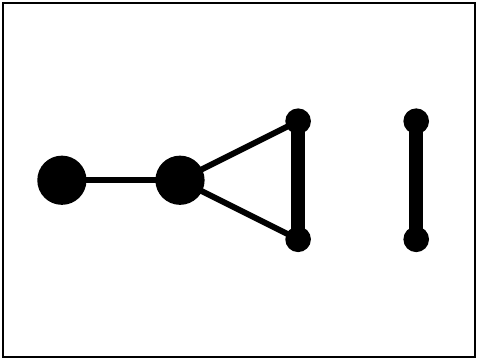}
\includegraphics[scale=0.25]{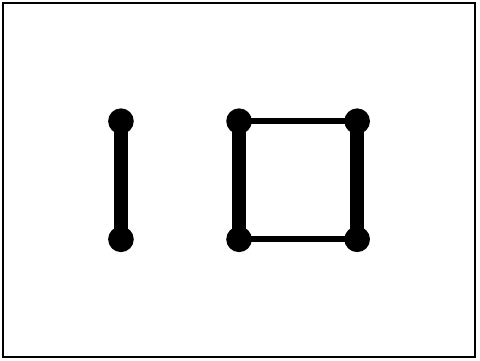}

\vspace{0.5ex}

\includegraphics[scale=0.25]{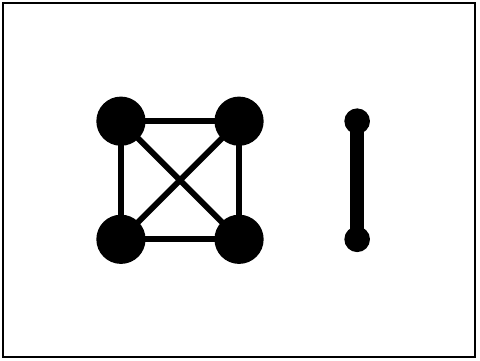}
\includegraphics[scale=0.25]{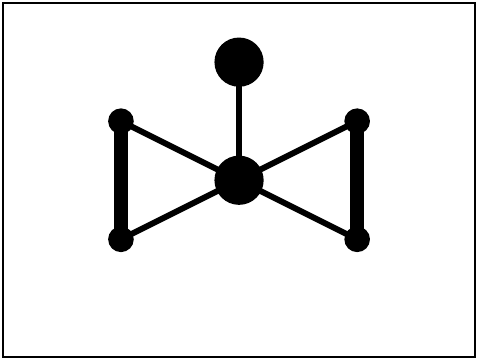}
\includegraphics[scale=0.25]{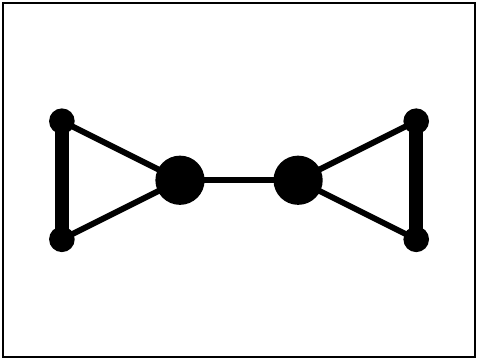}
\includegraphics[scale=0.25]{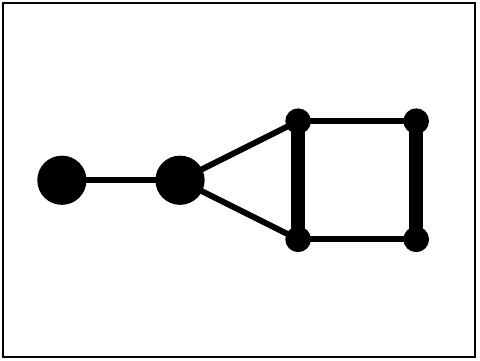}
\includegraphics[scale=0.25]{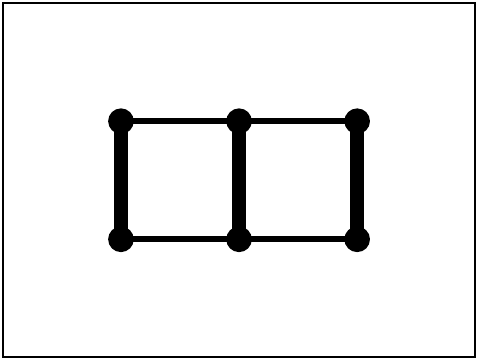}

\vspace{0.5ex}

\includegraphics[scale=0.25]{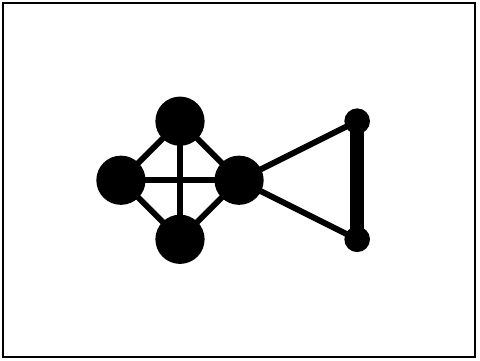}
\includegraphics[scale=0.25]{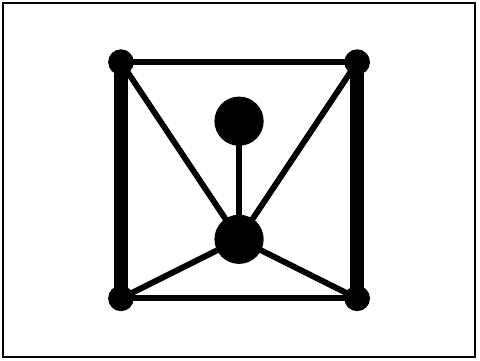}
\includegraphics[scale=0.25]{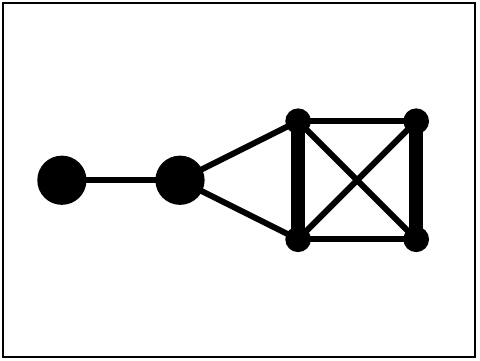}
\includegraphics[scale=0.25]{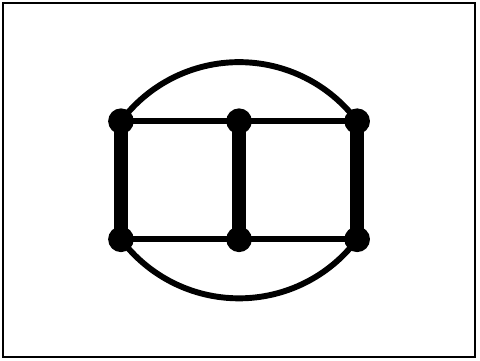}
\includegraphics[scale=0.25]{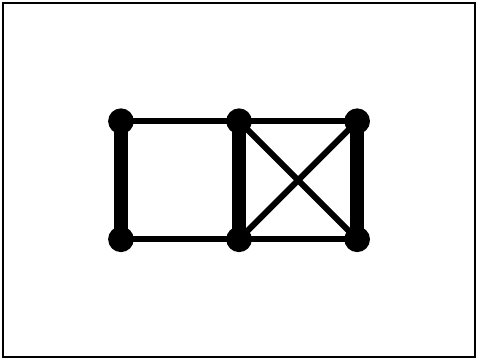}
\includegraphics[scale=0.25]{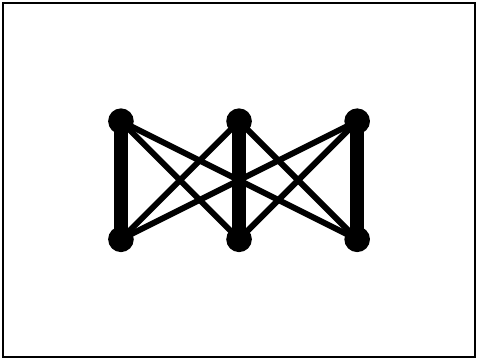}

\vspace{0.5ex}

\includegraphics[scale=0.25]{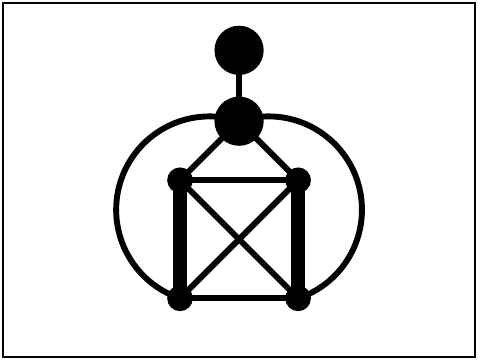}
\includegraphics[scale=0.25]{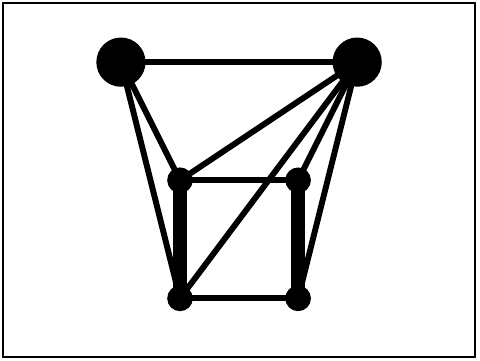}
\includegraphics[scale=0.25]{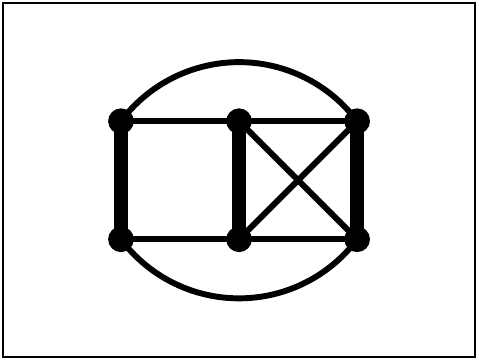}
\includegraphics[scale=0.25]{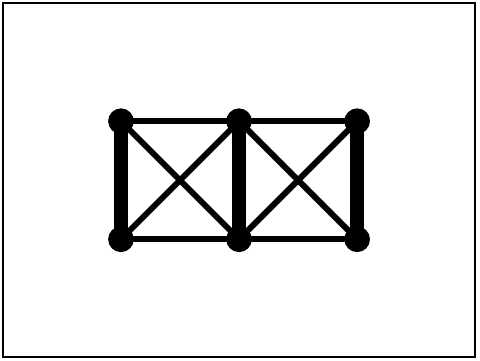}

\vspace{0.5ex}

\includegraphics[scale=0.25]{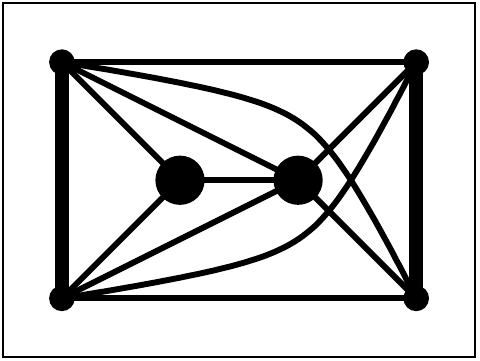}
\includegraphics[scale=0.25]{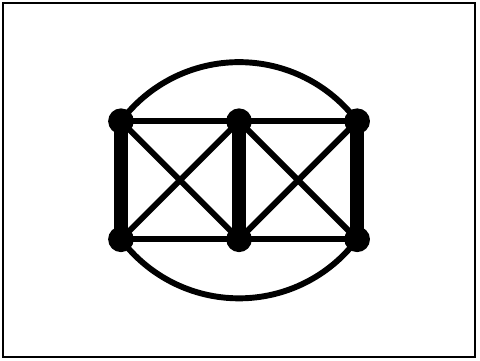}

\vspace{0.5ex}

\includegraphics[scale=0.25]{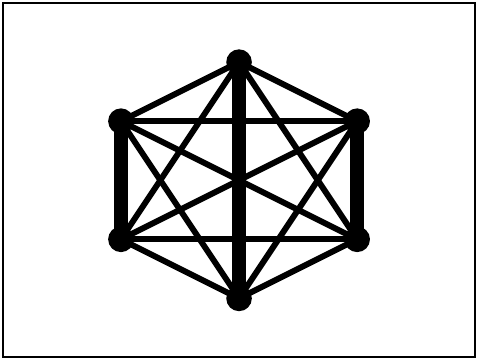}
\end{center}
\caption{\label{fig:sechs}The graphs admitting a strictly matched involution with 6 vertices,
ordered according to their numbers of edges: fat vertices indicate clique vertices and
fat edges indicate matching edges with respect to a certain (not necessarily unique) strictly 
matched involution}
\end{figure}

\cite{OCG} provided
a list of graphs admitting a strictly matched involution on $n$ vertices
with \begin{math}0\leq n\leq5\end{math}.
In Figure~\ref{fig:sechs}, as a further example, we display the graphs admitting a strictly
matched involution on 6 vertices.

\begin{table}[htbp]
\begin{center}
\scriptsize
\begin{tabular}{r|rrrrrrrrrrrr}
{\small$_m$ $^n$}&0&1&2&3&4&5&6&7&8&9&10&11\\
\hline
0&1&1&0&0&0&0&0&0&0&0&0&0\\
1&&&1&1&0&0&0&0&0&0&0&0\\
2&&&&0&1&1&0&0&0&0&0&0\\
3&&&&1&0&0&1&1&0&0&0&0\\
4&&&&&2&2&0&0&1&1&0&0\\
5&&&&&0&0&2&2&0&0&1&1\\
6&&&&&1&3&0&0&2&2&0&0\\
7&&&&&&0&5&5&0&0&2&2\\
8&&&&&&2&0&0&7&6&0&0\\
9&&&&&&0&6&11&0&0&7&6\\
10&&&&&&1&0&0&16&17&0&0\\
11&&&&&&&4&13&0&0&21&19\\
12&&&&&&&0&0&30&40&0&0\\
13&&&&&&&2&14&0&0&62&57\\
14&&&&&&&0&0&41&74&0&0\\
15&&&&&&&1&11&0&0&140&155\\
16&&&&&&&&0&50&128&0&0\\
17&&&&&&&&5&0&0&292&383\\
18&&&&&&&&0&40&175&0&0\\
19&&&&&&&&2&0&0&510&870\\
20&&&&&&&&0&28&199&0&0\\
21&&&&&&&&1&0&0&775&1745\\
22&&&&&&&&&16&184&0&0\\
23&&&&&&&&&0&0&987&3090\\
24&&&&&&&&&6&144&0&0\\
25&&&&&&&&&0&0&1098&4794\\
26&&&&&&&&&2&85&0&0\\
27&&&&&&&&&0&0&1009&6504\\
28&&&&&&&&&1&45&0&0\\
29&&&&&&&&&&0&801&7690\\
30&&&&&&&&&&19&0&0\\
31&&&&&&&&&&0&530&7956\\
32&&&&&&&&&&6&0&0\\
33&&&&&&&&&&0&302&7135\\
34&&&&&&&&&&2&0&0\\
35&&&&&&&&&&0&140&5596\\
36&&&&&&&&&&1&0&0\\
37&&&&&&&&&&&58&3791\\
38&&&&&&&&&&&0&0\\
39&&&&&&&&&&&20&2240\\
40&&&&&&&&&&&0&0\\
41&&&&&&&&&&&6&1138\\
42&&&&&&&&&&&0&0\\
43&&&&&&&&&&&2&507\\
44&&&&&&&&&&&0&0\\
45&&&&&&&&&&&1&195\\
46&&&&&&&&&&&&0\\
47&&&&&&&&&&&&67\\
48&&&&&&&&&&&&0\\
49&&&&&&&&&&&&21\\
50&&&&&&&&&&&&0\\
51&&&&&&&&&&&&6\\
52&&&&&&&&&&&&0\\
53&&&&&&&&&&&&2\\
54&&&&&&&&&&&&0\\
55&&&&&&&&&&&&1\\
\hline
$A(n)$&1&1&1&2&4&9&21&65&240&1128&6764&53971\\
\hline
\end{tabular}
\end{center}
\caption{\label{Anm}The number $C(n,m)$ of graphs admitting a strictly matched involution with $n$ vertices (and $m$ edges) for $n\leq11$}
\end{table}

\section{Counting by edges}\label{sec:exAB}

The proofs of
Observation~\ref{obs:A} and Theorem~\ref{obs:B}
are given in this section.
In order to simplify some arguments in the proofs,
we start with an obvious lemma.

\begin{lemma}\label{lem:hilfobs}
Every graph admitting a strictly matched involution with \begin{math}n\end{math} vertices
contains the graph \begin{math}\left\lfloor\frac{n}{2}\right\rfloor K_2\end{math} as a
(not necessarily induced) subgraph.
\end{lemma}

\begin{proof}
Let $G$ be a graph admitting a strictly matched involution with $n$ vertices.
By Theorem~\ref{thm:structure}, there is an integer $k$ with \begin{math}0\leq k\leq\nhalb\end{math}, such
that the vertex set of $G$ can be partitioned into a set containing a matching on $k$ edges
and a clique of size \begin{math}n-2k\end{math}. Since the clique contains a matching on 
\begin{math}\left\lfloor\frac{n-2k}{2}\right\rfloor\end{math} edges, the graph $G$ contains a
matching of size
\begin{displaymath}k+\left\lfloor\frac{n-2k}{2}\right\rfloor=\left\lfloor\frac{n}{2}\right\rfloor,\end{displaymath}
\textit{i.e.}, the graph \begin{math}\left\lfloor\frac{n}{2}\right\rfloor K_2\end{math} is a subgraph of~$G$.
\end{proof}

We will now prove Observation~\ref{obs:A} and Theorem~\ref{obs:B}.

\begin{proof}[of Observation~\ref{obs:A}]
We have the following.

{\bfseries Ad (\ref{oAa}):} 
Since, by Lemma~\ref{lem:hilfobs}, every graph admitting a 
strictly matched involution with $n$ vertices contains
the graph \begin{math}\left\lfloor\frac{n}{2}\right\rfloor K_2\end{math} as a subgraph,
which is a graph with \begin{math}\left\lfloor\frac{n}{2}\right\rfloor\end{math} edges,
(a)~follows trivially.

{\bfseries Ad (\ref{oAb}):} 
Since a graph on $n$ vertices has at most \begin{math}\binom{n}{2}\end{math} edges,
(b) follows trivially.

{\bfseries Ad (\ref{oAc}) and (\ref{oAd}):} 
Let $G$ be a graph admitting a strictly matched involution $\sigma$
with $n$ vertices, where \begin{math}n=4\ell+c\end{math} with \begin{math}\ell\in{\mathbb N}\end{math} and \begin{math}c\in\{0,1,2,3\}\end{math}.
Then the number \begin{math}m_0(n)\end{math} of matching edges and clique edges with respect to a fixed
strictly matched involution with $k$ matching edges is
\begin{align*}
m_0(n)&=k+\binom{n-2k}{2}\\
&=k+\frac{(4\ell+c-2k)(4\ell+c-2k-1)}{2}\\
&=k+\frac{4L}{2}+\frac{c^2-c+2k}{2}\\
&\equiv\binom{c}{2}\pmod{2},
\end{align*}
for some \begin{math}L\in{\mathbb Z}\end{math}. Thus, for \begin{math}c\in\{0,1\}\end{math},
\begin{displaymath}m_0(n)\equiv\binom{c}{2}\equiv0\pmod{2};\end{displaymath}
and, for \begin{math}c\in\{2,3\}\end{math},
\begin{displaymath}m_0(n)\equiv\binom{c}{2}\equiv1\pmod{2}.\end{displaymath}
Since all other types of edges are added pairwisely, the number of
edges of $G$ is even if \begin{math}c\in\{0,1\}\end{math}, \textit{i.e.}, if \begin{math}n\equiv0,1\end{math} (mod 4);
and the number of edges of \begin{math}G\end{math} is odd if \begin{math}c\in\{2,3\}\end{math}, \textit{i.e.}, if \begin{math}n\equiv2,3\end{math} (mod 4). This proves (\ref{oAc}) and (\ref{oAd}). 
\end{proof}

\pagebreak[3]

\begin{proof}[of Theorem~\ref{obs:B}]
We have the following.

{\bfseries Ad (\ref{oBa}):} 
Let $G$ be a graph admitting a strictly matched involution with
$n$ vertices and \begin{math}\nhalb\end{math} edges. 
By Lemma~\ref{lem:hilfobs}, $G$~has the graph \begin{math}\nhalb K_2\end{math} as a subgraph, which is a graph
with \begin{math}\nhalb\end{math} edges.
Thus, 
$G$ is isomorphic to \begin{math}\nhalb K_2\end{math} if $n$ is even, respectively,
$G$ is isomorphic to \begin{math}\nhalb K_2\cup K_1\end{math} if $n$ is odd.
Both types of graphs are obviously 
graphs admitting a strictly matched involution (with clique size 0,
respectively, clique size 1).
This proves~(\ref{oBa}).

{\bfseries Ad (\ref{oBb}):} 
Let $G$ be a graph admitting a strictly matched involution
with $n$ vertices and \begin{math}m=\nhalb+2\end{math} edges.
By Lemma~\ref{lem:hilfobs}, $G$ contains \begin{math}\nhalb K_2\end{math} as a subgraph.
Let $G'$ be the subgraph on the vertex set of $G$ with the edges of the
\begin{math}\nhalb K_2\end{math}.
Furthermore there are two additional edges in $G$. 

If one of these edges is adjacent to matching edges $e$ and~$f$, then, 
by Theorem~\ref{thm:structure}, the other edge is also adjacent to $e$ and~$f$,
moreover if $n$ is odd, then the other edge must be symmetrical to the first edge.
If one of these edges is adjacent to a matching edge $e$ and 
incident with the isolated vertex of~$G'$ 
(in case $n$ is odd), then, by Theorem~\ref{thm:structure}, the other edge
must also be adjacent to $e$ and incident with the isolated vertex of~$G'$.

We conclude that in case $n$ is even (and \begin{math}n\geq4\end{math}), the graph $G$ is
isomorphic either to 
\begin{math}C_4\cup\left(\nhalb-2\right)K_2\end{math}
or to
\begin{math}P\cup\left(\nhalb-2\right)K_2\end{math},
where $P$ is the paw graph.
And in case $n$ is odd (and \begin{math}n\geq5\end{math}), the graph $G$ 
is isomorphic either to
\begin{math}K_1\cup C_4\cup\left(\nhalb-2\right)K_2\end{math}
or to
\begin{math}K_3\cup\left(\nhalb-1\right)K_2\end{math}.

As all these graphs are obviously graphs admitting a strictly matched involution,
(\ref{oBb})~is proved.
 
{\bfseries Ad (\ref{oBc}):} 
Let $n$ be an even number, \begin{math}n\geq8\end{math}, and
$G$ be a graph admitting a strictly matched involution
with $n$ vertices and \begin{math}m=\frac{n}{2}+4\end{math} edges.
By Theorem~\ref{thm:structure}, there is \begin{math}k\in{\mathbb{N}}\end{math}, so that
the vertex set of $G$ can be decomposed into a clique $C$ of size \begin{math}s=n-2k\end{math}
and the vertex sets of two isomorphic graphs $H$ and $H'$, 
each with $k$ vertices, moreover the clique size is even and cannot be
more than~4, thus for $G$ only the following cases are possible:
\begin{enumerate}[(1)]
\item \label{oBcA}\begin{math}s=0\end{math} and \begin{math}k=\frac{n}{2}\end{math}
\item \label{oBcB}\begin{math}s=2\end{math} and \begin{math}k=\frac{n}{2}-1\end{math}
\item \label{oBcC}\begin{math}s=4\end{math} and \begin{math}k=\frac{n}{2}-2\end{math}
\end{enumerate}
If we let \begin{math}m(H)\end{math} denote the number of edges of $H$, \begin{math}m(C)\end{math} denote the
number of edges of the graph induced by the clique $C$, 
\begin{math}m_a\end{math} denote the number of pairs of edges between $H$ and $H'$
which are not matching edges (we call them \emph{cross-edges}), 
and \begin{math}m_b\end{math} denote the number of pairs of edges between $H$ and $C$
(we call them \emph{clique-connecting-edges}),
then we note that in Case (\ref{oBcA}) and in Case (\ref{oBcB})
\begin{displaymath}m(H)+m_a+m_b=2,\end{displaymath}
whereas in Case (\ref{oBcC})
\begin{displaymath}m(H)+m_a+m_b=0.\end{displaymath}
Therefore,
we subdivide the cases (\ref{oBcA}), (\ref{oBcB}), and (\ref{oBcC}) as follows.
\begin{enumerate}[(\ref{oBcA}.a)]
\item \label{oBcAa}\begin{math}m(C)=0\end{math} and \begin{math}m(H)=0\end{math} (thus \begin{math}m_a=2\end{math} and \begin{math}m_b=0\end{math})
\item \label{oBcAb}\begin{math}m(C)=0\end{math} and \begin{math}m(H)=1\end{math} (thus \begin{math}m_a=1\end{math} and \begin{math}m_b=0\end{math})
\item \label{oBcAc}\begin{math}m(C)=0\end{math} and \begin{math}m(H)=2\end{math} (thus \begin{math}m_a=0\end{math} and \begin{math}m_b=0\end{math})
\end{enumerate}

\begin{enumerate}[(\ref{oBcB}.a)]
\item \label{oBcBa}\begin{math}m(C)=1\end{math} and \begin{math}m(H)=0\end{math} and \begin{math}m_a=2\end{math} (thus \begin{math}m_b=0\end{math})
\item \label{oBcBb}\begin{math}m(C)=1\end{math} and \begin{math}m(H)=0\end{math} and \begin{math}m_a=1\end{math} (thus \begin{math}m_b=1\end{math})
\item \label{oBcBc}\begin{math}m(C)=1\end{math} and \begin{math}m(H)=0\end{math} and \begin{math}m_a=0\end{math} (thus \begin{math}m_b=2\end{math})
\item \label{oBcBd}\begin{math}m(C)=1\end{math} and \begin{math}m(H)=1\end{math} and \begin{math}m_a=1\end{math} (thus \begin{math}m_b=0\end{math})
\item \label{oBcBe}\begin{math}m(C)=1\end{math} and \begin{math}m(H)=1\end{math} and \begin{math}m_a=0\end{math} (thus \begin{math}m_b=1\end{math})
\item \label{oBcBf}\begin{math}m(C)=1\end{math} and \begin{math}m(H)=2\end{math} (thus \begin{math}m_a=0\end{math} and \begin{math}m_b=0\end{math})
\end{enumerate}

\begin{enumerate}[(\ref{oBcC}.a)]
\item \label{oBcCa}\begin{math}m(C)=6\end{math} (thus \begin{math}m(H)=0\end{math} and \begin{math}m_a=0\end{math} and \begin{math}m_b=0\end{math})
\end{enumerate}

In the following, we refer to the names of types of graphs depicted
in Figure~\ref{fig:graphH}.

\begin{figure}[tbp]
\begin{center}
\begin{minipage}[t]{0.36\textwidth}
\begin{math}G_1\end{math}\quad\includegraphics[scale=0.3]{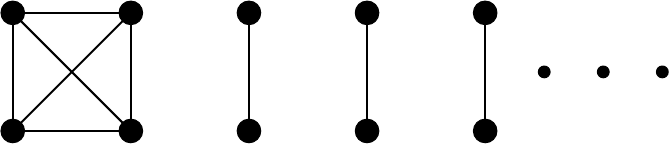}

\bigskip

\begin{math}G_2\end{math}\quad\includegraphics[scale=0.3]{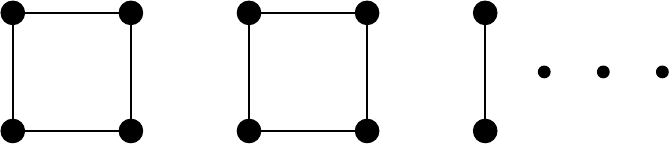}

\bigskip

\begin{math}G_3\end{math}\quad\includegraphics[scale=0.3]{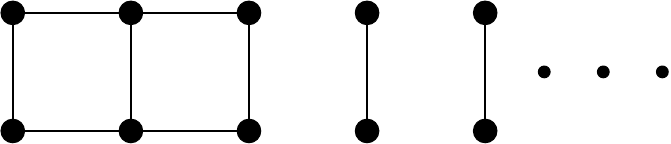}

\bigskip

\begin{math}G_4\end{math}\quad\includegraphics[scale=0.3]{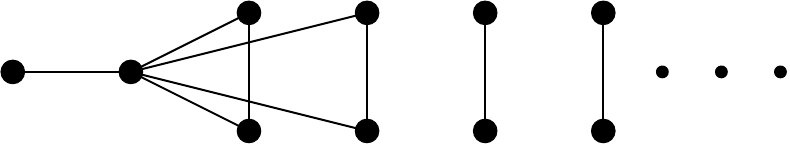}

\bigskip

\begin{math}G_5\end{math}\quad\includegraphics[scale=0.3]{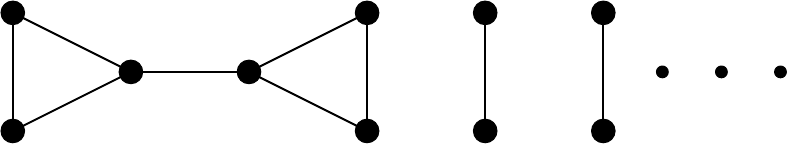}

\bigskip

\begin{math}G_6\end{math}\quad\includegraphics[scale=0.3]{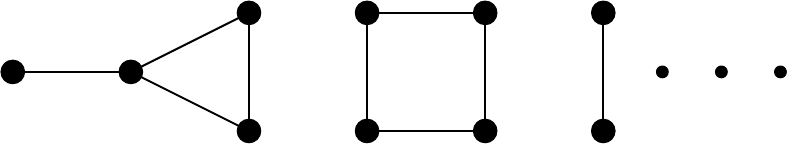}

\bigskip

\begin{math}G_7\end{math}\quad\includegraphics[scale=0.3]{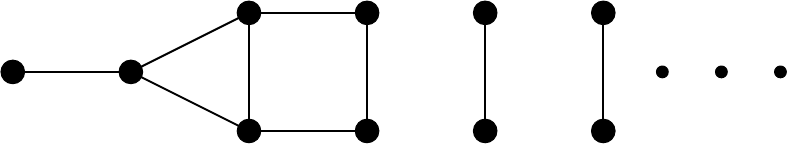}
\end{minipage}\qquad\quad
\begin{minipage}[t]{0.36\textwidth}
\begin{math}I_1\end{math}\quad\includegraphics[scale=0.3]{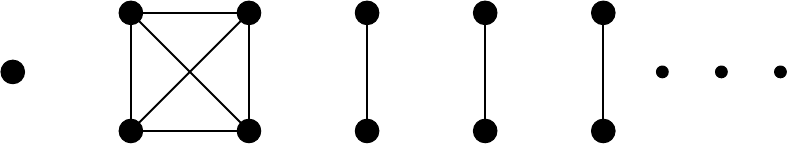}

\bigskip

\begin{math}I_2\end{math}\quad\includegraphics[scale=0.3]{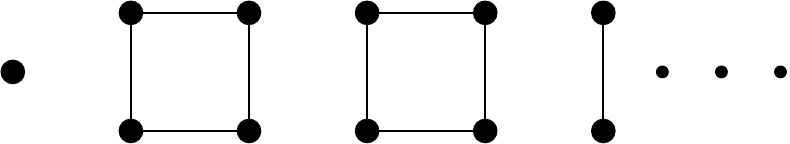}

\bigskip

\begin{math}I_3\end{math}\quad\includegraphics[scale=0.3]{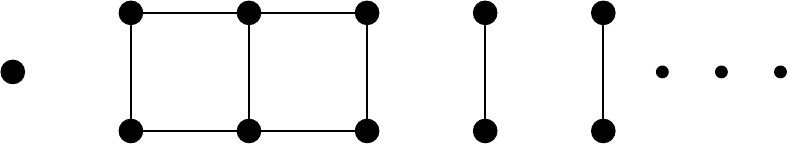}

\bigskip

\begin{math}I_4\end{math}\quad\includegraphics[scale=0.3]{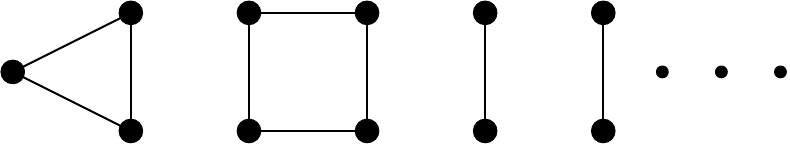}

\bigskip

\begin{math}I_5\end{math}\quad\includegraphics[scale=0.3]{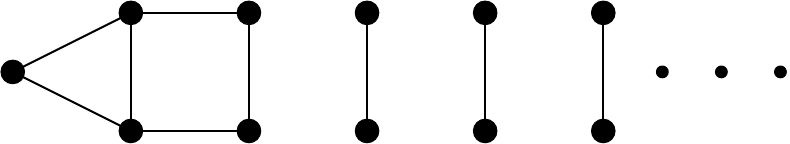}

\bigskip

\begin{math}I_6\end{math}\quad\includegraphics[scale=0.3]{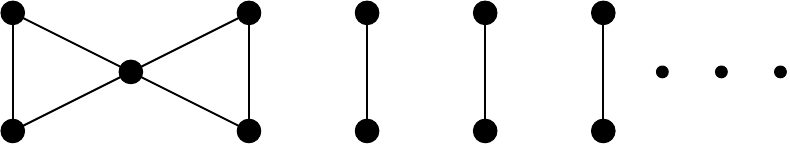}
\end{minipage}
\end{center}
\caption{\label{fig:graphH}The graphs occurring in the proof of Theorem~\ref{obs:B}~(c) and~(d)}
\end{figure}

In {\bf Case (\ref{oBcA}.\ref{oBcAa})}, depending on whether all four cross-edges are adjacent
to one certain matching edge or not, we see that \begin{math}G\end{math} is isomorphic to $G_3$ or $G_2$,
respectively.

In {\bf Case (\ref{oBcA}.\ref{oBcAb})}, 
if both cross-edges are adjacent to the edge of $H$,
the graph $G$ is isomorphic to $G_1$. If only one cross-edge is adjacent
to the edge of $H$, then $G$ is isomorphic to $G_3$. If none of the cross-edges
is adjacent to the edge of $H$, then $G$ is isomorphic to $G_2$.

In {\bf Case (\ref{oBcA}.\ref{oBcAc})}, 
if the two edges of $H$ are connected, we obtain $G_3$,
otherwise we obtain $G_2$.

In {\bf Case (\ref{oBcB}.\ref{oBcBa})}, 
we obtain the same as in Case (\ref{oBcA}.\ref{oBcAa}), since the clique edge
cannot be connected by cross-edges, thus the clique edge can be regarded
as a matching edge here.

In {\bf Case (\ref{oBcB}.\ref{oBcBb})}, 
let $e$ be the matching edge that is adjacent to
the two clique-connecting-edges. If the two cross-edges are adjacent
to $e$, then we obtain $G_7$, otherwise we obtain $G_6$ as our graph~$G$.

In {\bf Case (\ref{oBcB}.\ref{oBcBc})}, 
if the four clique-connecting-edges are incident
with the same vertex of the clique, we obtain $G_4$, otherwise we get~$G_5$.

In {\bf Case (\ref{oBcB}.\ref{oBcBd})}, 
we obtain mainly the same as in Case (\ref{oBcA}.\ref{oBcAb}), since the
clique edge can be regarded again as a matching edge (but note that $G_2$
can appear in Case (\ref{oBcB}.\ref{oBcBd}) only if $n\geq10$).

In {\bf Case (\ref{oBcB}.\ref{oBcBe})}, 
if the edge of $H$ is adjacent to a clique-connecting-edge,
then we obtain $G_7$. Otherwise, we obtain~$G_6$.

In {\bf Case (\ref{oBcB}.\ref{oBcBf})}, 
we obtain the same as in Case (\ref{oBcA}.\ref{oBcAc}), since the
clique edge can be regarded again as a matching edge.

In {\bf Case (\ref{oBcC}.\ref{oBcCa})}, we obtain $G_1$.

Thus, on the one hand, 
each of the types of graphs on the left-hand side of Figure~\ref{fig:graphH}
may appear. On the other hand, these types of graphs admit a strictly matched
involution (\textit{e.g.}, with the matching edges drawn vertically in the figure).

{\bfseries Ad (\ref{oBd}):} 
Let $n$ be an odd number, \begin{math}n\geq9\end{math}, and
$G$ be a graph admitting a strictly matched involution
with $n$ vertices and \begin{math}m=\frac{n-1}{2}+4\end{math} edges.
By Theorem~\ref{thm:structure}, there is \begin{math}k\in{\mathbb{N}}\end{math}, so that
the vertex set of $G$ can be decomposed into a clique $C$ of size \begin{math}s=n-2k\end{math}
and the vertex sets of two isomorphic graphs $H$ and $H'$, 
each with $k$ vertices, moreover the clique size is odd
and cannot be more than~4, thus for $G$ only the following cases are possible:
\begin{enumerate}[(1)]
\item \label{oBdA}\begin{math}s=1\end{math} and \begin{math}k=\frac{n-1}{2}\end{math}
\item \label{oBdB}\begin{math}s=3\end{math} and \begin{math}k=\frac{n-1}{2}-1\end{math}
\end{enumerate}
With the same notation as in the proof of (\ref{oBc}),
we note that in Case (\ref{oBdA})
\begin{displaymath}m(H)+m_a+m_b=2,\end{displaymath}
whereas in Case (\ref{oBdB})
\begin{displaymath}m(H)+m_a+m_b=1.\end{displaymath}
Thus,
we subdivide the cases (\ref{oBdA}) and (\ref{oBdB}) as follows.
\begin{enumerate}[(\ref{oBdA}.a)]
\item \label{oBdAa}\begin{math}m(C)=0\end{math} and \begin{math}m(H)=0\end{math} and \begin{math}m_a=2\end{math} (thus \begin{math}m_b=0\end{math})
\item \label{oBdAb}\begin{math}m(C)=0\end{math} and \begin{math}m(H)=0\end{math} and \begin{math}m_a=1\end{math} (thus \begin{math}m_b=1\end{math})
\item \label{oBdAc}\begin{math}m(C)=0\end{math} and \begin{math}m(H)=0\end{math} and \begin{math}m_a=0\end{math} (thus \begin{math}m_b=2\end{math})
\item \label{oBdAd}\begin{math}m(C)=0\end{math} and \begin{math}m(H)=1\end{math} and \begin{math}m_a=1\end{math} (thus \begin{math}m_b=0\end{math})
\item \label{oBdAe}\begin{math}m(C)=0\end{math} and \begin{math}m(H)=1\end{math} and \begin{math}m_a=0\end{math} (thus \begin{math}m_b=1\end{math})
\item \label{oBdAf}\begin{math}m(C)=0\end{math} and \begin{math}m(H)=2\end{math} (thus \begin{math}m_a=0\end{math} and \begin{math}m_b=0\end{math})
\end{enumerate}

\begin{enumerate}[(\ref{oBdB}.a)]
\item \label{oBdBa}\begin{math}m(C)=3\end{math} and \begin{math}m(H)=0\end{math} and \begin{math}m_a=1\end{math} (thus \begin{math}m_b=0\end{math})
\item \label{oBdBb}\begin{math}m(C)=3\end{math} and \begin{math}m(H)=0\end{math} and \begin{math}m_a=0\end{math} (thus \begin{math}m_b=1\end{math})
\item \label{oBdBc}\begin{math}m(C)=3\end{math} and \begin{math}m(H)=1\end{math} (thus \begin{math}m_a=0\end{math} and \begin{math}m_b=0\end{math})
\end{enumerate}

\pagebreak[3]

In the following, we refer again to the names of types of graphs depicted
in Figure~\ref{fig:graphH}.

In {\bf Case (\ref{oBdA}.\ref{oBdAa})}, as \begin{math}m_b=0\end{math}, by the same arguments 
as in Case (\ref{oBcA}.\ref{oBcAa}) in the proof of~(\ref{oBc}),
we obtain $I_2$ or $I_3$ as graph~$G$.

In {\bf Case (\ref{oBdA}.\ref{oBdAb})}, if a cross-edge is adjacent to a
clique-connecting edge, then we obtain~$I_5$. Otherwise,
we obtain~$I_4$.

In {\bf Case (\ref{oBdA}.\ref{oBdAc})} we necessarily obtain $I_6$.

In {\bf Case (\ref{oBdA}.\ref{oBdAc})}, as \begin{math}m_b=0\end{math}, by the same arguments 
as in Case (\ref{oBcA}.\ref{oBcAb}) in the proof of~(\ref{oBc}),
we obtain $I_1$ or $I_2$ or $I_3$.

In {\bf Case (\ref{oBdA}.\ref{oBdAe})}, if the edge of $H$ is adjacent to a
clique-connecting-edge, then we obtain~$I_5$. Otherwise,
we obtain~$I_4$.

In {\bf Case (\ref{oBdA}.\ref{oBdAf})}, as \begin{math}m_b=0\end{math}, by the same arguments 
as in Case (\ref{oBcA}.\ref{oBcAc}) in the proof of~(\ref{oBc}),
we obtain $I_2$ or $I_3$.

In {\bf Case (\ref{oBdB}.\ref{oBdBa})} we obtain $I_4$.

In {\bf Case (\ref{oBdB}.\ref{oBdBb})} we obtain $I_6$.

In {\bf Case (\ref{oBdB}.\ref{oBdBc})} we obtain $I_4$ again.

Thus, on the one hand, 
each of the types of graphs on the right-hand side of Figure~\ref{fig:graphH}
may appear. On the other hand, these types of graphs admit a strictly matched
involution (\textit{e.g.}, with the matching edges drawn vertically in the figure). 
\end{proof}

The proofs of Theorem~\ref{obs:B} (\ref{oBc}) and (\ref{oBd}) may give a good idea for the
fact that in our upper bound of Theorem~\ref{oberneu} and in the enumeration
approach derived from Theorem~\ref{oberneu} many graphs are counted multiple times, which is getting
super-exponentially worse with increasing number of vertices and edges. This
is the reason why in our computer program isomorphism tests consume most
of the running time.

Using Observation~\ref{obs:A} and Theorem~\ref{obs:B} and the values from Table~\ref{tab:an},
we can determine \begin{math}B(m)\end{math} for \begin{math}0\leq m\leq11\end{math}, see Table~\ref{AnmB}.
We have implemented these values \begin{math}B(m)\end{math} in the OEIS (\cite{oeis}) in
sequence \seqnum{A363772}.

\begin{table}[htbp]
\begin{center}
\begin{tabular}{r|rrrrrrrrrrrr}
{\small$_n$ $^m$}&0&1&2&3&4&5&6&7&8&9&10&11\\
\hline
0&1&&&&&&&&&&&\\
1&1&&&&&&&&&&&\\
2&&1&&&&&&&&&&\\
3&&1&0&1&&&&&&&&\\
4&&&1&0&2&0&1&&&&&\\
5&&&1&0&2&0&3&0&2&0&1&\\
6&&&&1&0&2&0&5&0&6&0&4\\
7&&&&1&0&2&0&5&0&11&0&13\\
8&&&&&1&0&2&0&7&0&16&0\\
9&&&&&1&0&2&0&6&0&17&0\\
10&&&&&&1&0&2&0&7&0&21\\
11&&&&&&1&0&2&0&6&0&19\\
12&&&&&&&1&0&2&0&7&0\\
13&&&&&&&1&0&2&0&6&0\\
14&&&&&&&&1&0&2&0&7\\
15&&&&&&&&1&0&2&0&6\\
16&&&&&&&&&1&0&2&0\\
17&&&&&&&&&1&0&2&0\\
18&&&&&&&&&&1&0&2\\
19&&&&&&&&&&1&0&2\\
20&&&&&&&&&&&1&0\\
21&&&&&&&&&&&1&0\\
22&&&&&&&&&&&&1\\
23&&&&&&&&&&&&1\\
\hline
$B(m)$&2&2&2&3&6&6&10&16&21&36&53&76\\
\hline
\end{tabular}
\end{center}
\caption{\label{AnmB}The number of graphs admitting a strictly matched involution with $m$ edges (and $n$ vertices) for $m\leq11$.
The values for $0\leq n\leq11$ are taken from Table~\ref{Anm}, the values for
$12\leq n\leq 23$ are implied by Observation~\ref{obs:A} and Theorem~\ref{obs:B}.}
\end{table}

\section{The decomposition lemma}\label{sec:exC}

The key fact for the proof of Theorem~\ref{obs:C} will be the
following decomposition lemma.
We start with three formulations of this main lemma. Its simple form
is the following.

\begin{lemma}[decomposition lemma, version 1: simple form]\label{lem:decA}
For every graph $G$, the following statements are equivalent:
\begin{enumerate}[(1)]
\item $G$ admits a strictly matched involution.
\item \begin{math}G\vee K_1\end{math} admits a strictly matched involution.
\end{enumerate}
\end{lemma}

By induction, we get the second form.

\begin{lemma}[decomposition lemma, version 2: aggregated form]\label{lem:decB}
For every graph $G$ and every \begin{math}n\in{\mathbb N}\end{math}, the following statements are equivalent:
\begin{enumerate}[(1)]
\item \label{decLemA}$G$ admits a strictly matched involution.
\item \label{decLemB}\begin{math}G\vee K_n\end{math} admits a strictly matched involution.
\end{enumerate}
\end{lemma}

We will use a third form of the decomposition lemma, which is
derived from the second form.

\begin{lemma}[decomposition lemma, version 3: complementary form]\label{lem:decC}
For every $n\in{\mathbb N}$ and every subgraph \begin{math}S=(V,E)\end{math} of \begin{math}K_n\end{math}, 
the following statements are equivalent:
\begin{enumerate}[(1)]
\item The complement \begin{math}\overline{S}\end{math} of $S$ admits a strictly matched involution.
\item \begin{math}K_n-E\end{math} admits a strictly matched involution.
\end{enumerate}
\end{lemma}

We will now prove the simple form of the decomposition lemma, \textit{i.e.}, Lemma~\ref{lem:decA}.
The other forms (Lemmata~\ref{lem:decB} and~\ref{lem:decC}) can be obtained as easy
corollaries from Lemma~\ref{lem:decA}.

\begin{proof}[of Lemma~\ref{lem:decA}]
Assume $G$ admits a strictly matched involution $\sigma$. If we add a new fixed point $v$
to $G$ that is connected to every vertex of $G$, then the mapping $\sigma'$ with
\begin{displaymath}
\sigma'(w) = \begin{cases}
    \sigma(w), & \text{if \begin{math}w\neq v\end{math};}\\
    v,         & \text{if \begin{math}w=v\end{math};}
\end{cases}
\end{displaymath}
is a strictly matched, involutive graph automorphism. This proves the trivial implication
\begin{math}(\ref{decLemA})\Longrightarrow(\ref{decLemB})\end{math}.

Now assume that \begin{math}G'=G\vee K_1\end{math} admits a strictly matched involution \begin{math}\sigma'\end{math}. Let $v$ be the vertex
of $K_1$. If $v$ is in the fixed point set \begin{math}F_1(G',\sigma')\end{math}, then we may obviously delete $v$
and \begin{math}\sigma_{|G}'\end{math} is still strictly matched, well-defined, and involutive. Thus, suppose that
\begin{displaymath}v\in F_2(G',\sigma').\end{displaymath}
Let \begin{math}w=\sigma'(v)\end{math}. Since $v$ is connected to every vertex of $G$, 
the vertex $w$ is connected to every other vertex of $G'$ 
(as \begin{math}\sigma'\end{math} is a graph automorphism).
Thus we may
consider $w$ as a part of the clique, \textit{i.e.}, we may define a map $\sigma$ on the vertex set of
$G$ by
\begin{displaymath}
\sigma(x) = \begin{cases}
    \sigma'(x), & \text{if \begin{math}x\neq w\end{math};}\\
    w,          & \text{if \begin{math}x=w\end{math}.}
\end{cases}
\end{displaymath}
The map $\sigma$ is a strictly matched involutive graph automorphism, since \begin{math}\sigma'\end{math} has the same
property on \begin{math}G-w\end{math} and $w$ is a fixed point of $\sigma$ connected to every other vertex, in particular,
to every pair of vertices of orbit~2. This proves the non-trivial implication 
\begin{math}(\ref{decLemB})\Longrightarrow(\ref{decLemA})\end{math}.
\end{proof}

\begin{proof}[of Lemma~\ref{lem:decB}]
Applying Lemma~\ref{lem:decA} inductively $n$ times implies the statement
of Lemma~\ref{lem:decB}.
\end{proof}

\begin{proof}[of Lemma~\ref{lem:decC}]
Since 
\begin{math}K_n-E\end{math} is isomorphic to \begin{math}\overline{S}\vee K_{n-|V|}\end{math},
the statement of Lemma~\ref{lem:decC} directly follows from Lemma~\ref{lem:decB}
when substituting $G$ by \begin{math}\overline{S}\end{math} and $n$ by \begin{math}n-|V|\end{math}.
\end{proof}

In the proof of Theorem~\ref{obs:C} we will also frequently make
use of the following obvious lemma.

\begin{lemma}\label{lem:degree}
For every graph \begin{math}G\end{math} admitting a strictly matched involution $\sigma$ and
every edge \begin{math}vw\end{math} of $G$ with \begin{math}\sigma(v)=w\end{math}, \textit{i.e.}, for every matching edge~\begin{math}vw\end{math}, we have
\begin{displaymath}d_G(v)=d_G(w).\end{displaymath}
\end{lemma}

Using Lemma~\ref{lem:decC} and
Lemma~\ref{lem:degree}
we are able to give the following simple
proof of Theorem~\ref{obs:C}.

\begin{proof}[of Theorem~\ref{obs:C}]
We have the following.

{\bfseries Ad (\ref{oCa}):}
Since the only graph with $n$ vertices and \begin{math}\binom{n}{2}\end{math} 
edges is the complete graph on $n$ vertices, which is a graph admitting a
strictly matched involution with \begin{math}k=0\end{math} matching edges and a clique of
size $n$, the assertion follows trivially.

{\bfseries Ad (\ref{oCb}):}
For \begin{math}n\geq4\end{math}, there are only two graphs with \begin{math}\binom{n}{2}-2\end{math} edges.
They can be obtained from the complete graph on $n$ edges
by deleting either 
two nonadjacent or two adjacent edges.
We will prove that these two graphs admit a strictly matched involution.

If $n$ is even, then the first graph can be constructed by using a clique of
size 0 and adding all possible edges except for one pair of edges between
two matching edges; and the second graph can be constructed by using a clique
of size 2 and adding all possible edges except for one pair of edges between
the two vertices of a matching edge and one of the vertices of the clique. 

Similarly, if $n$ is odd, then the first graph can be constructed by using a clique of
size 1 and adding all possible edges except for one pair of edges between
two matching edges; and the second graph can be constructed by using a clique
of size 1 and adding all possible edges except for one pair of edges between
the two vertices of a matching edge and the vertex of the clique. 

\begin{figure}[htbp]
\begin{center}\small
\begin{tabular}{|lcr|lcr|c|}
\hline
&&&&&&\\[-1.6ex]
&&&&&&$\overline{S}$ admits a strictly\\
&$S$&&&$\overline{S}$&&matched involution?\\
&&&&&&\\[-2.1ex]
\hline\hline
&&&&&&\\[-1.6ex]
$S_1$&\includegraphics[scale=0.3]{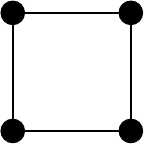}&
&$\overline{S_1}$&\includegraphics[scale=0.3]{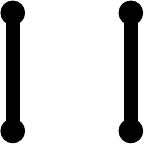}&
&YES\\
&&&&&&\\[-2.1ex]
\hline
&&&&&&\\[-1.6ex]
$S_2$&\includegraphics[scale=0.3]{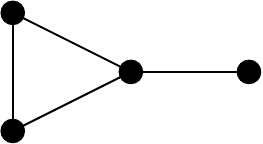}&
&$\overline{S_2}$&\includegraphics[scale=0.3]{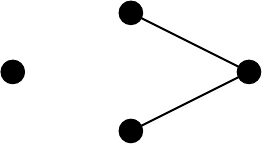}&
&NO\\
&&&&&&\\[-2.1ex]
\hline
&&&&&&\\[-1.6ex]
$S_3$&\includegraphics[scale=0.3]{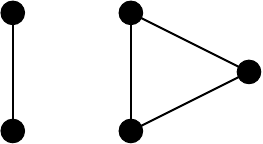}&
&$\overline{S_3}$&\includegraphics[scale=0.3]{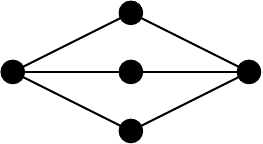}&
&NO\\
&&&&&&\\[-2.1ex]
\hline
&&&&&&\\[-1.6ex]
$S_4$&\includegraphics[scale=0.3]{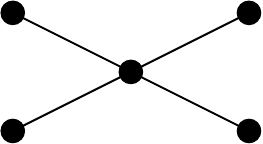}&
&$\overline{S_4}$&\includegraphics[scale=0.3]{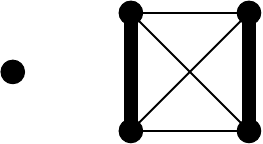}&
&YES\\
&&&&&&\\[-2.1ex]
\hline
&&&&&&\\[-1.6ex]
$S_5$&\includegraphics[scale=0.3]{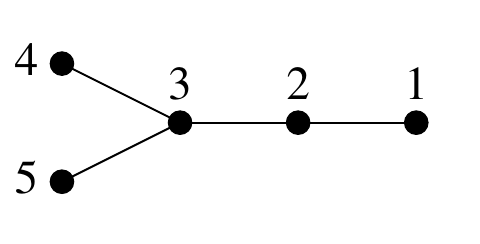}&
&$\overline{S_5}$&\includegraphics[scale=0.3]{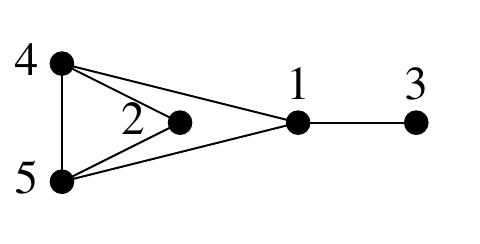}&
&NO\\
&&&&&&\\[-2.1ex]
\hline
&&&&&&\\[-1.6ex]
$S_6$&\includegraphics[scale=0.3]{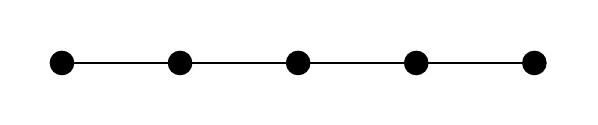}&
&$\overline{S_6}$&\includegraphics[scale=0.3]{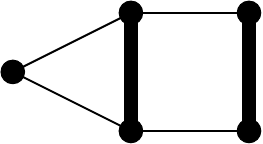}&
&YES\\
&&&&&&\\[-2.1ex]
\hline
&&&&&&\\[-1.6ex]
$S_7$&\includegraphics[scale=0.3]{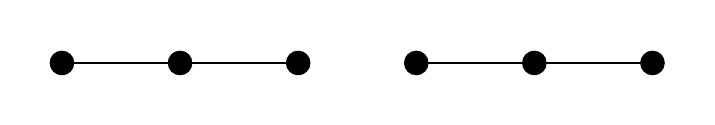}&
&$\overline{S_7}$&\includegraphics[scale=0.3]{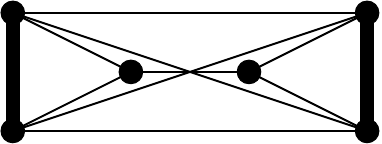}&
&YES\\
&&&&&&\\[-2.1ex]
\hline
&&&&&&\\[-1.6ex]
$S_8$&\includegraphics[scale=0.3]{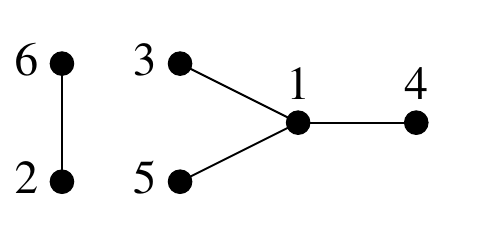}&
&$\overline{S_8}$&\includegraphics[scale=0.3]{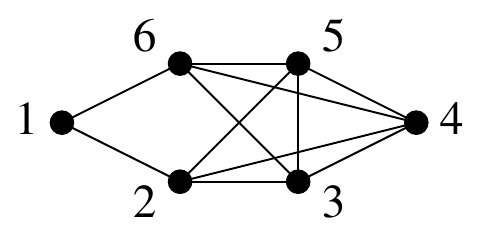}&
&NO\\
&&&&&&\\[-2.1ex]
\hline
&&&&&&\\[-1.6ex]
$S_9$&\includegraphics[scale=0.3]{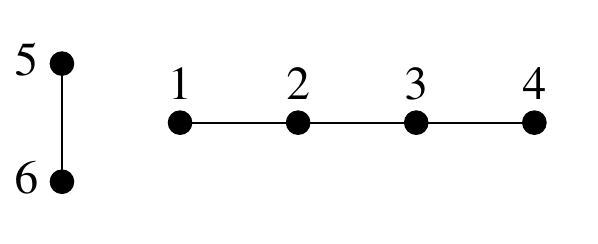}&
&$\overline{S_9}$&\includegraphics[scale=0.3]{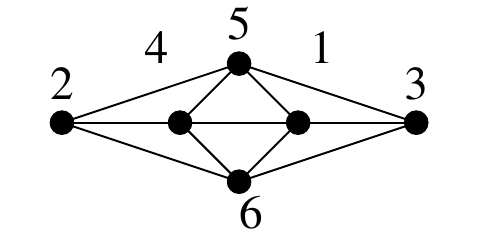}&
&NO\\
&&&&&&\\[-2.1ex]
\hline
&&&&&&\\[-1.6ex]
$S_{10}$&\includegraphics[scale=0.3]{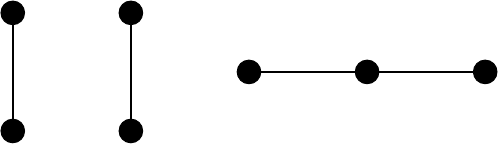}&
&$\overline{S_{10}}$&\includegraphics[scale=0.3]{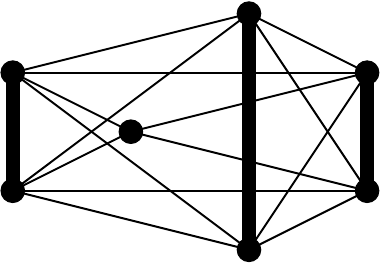}&
&YES\\
&&&&&&\\[-2.1ex]
\hline
&&&&&&\\[-1.6ex]
$S_{11}$&\includegraphics[scale=0.3]{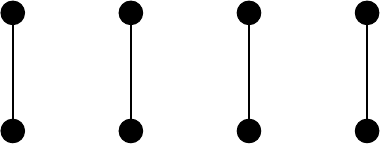}&
&$\overline{S_{11}}$&\includegraphics[scale=0.3]{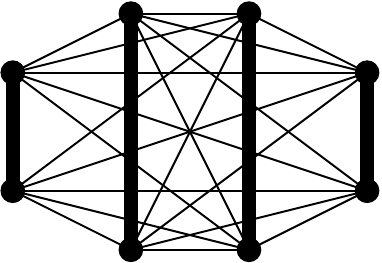}&
&YES\\[-2.1ex]
&&&&&&\\
\hline
\end{tabular}
\end{center}
\caption{\label{fig:minusfour}List of all graphs $S_i$ with 4 edges 
(without isolated vertices) and their complements $\overline{S_i}$}
\end{figure}

{\bfseries Ad (\ref{oCc}):}
Let $G$ be a graph admitting a strictly matched involution
with $n$ vertices, where $n\geq8$, and 
with \begin{math}m=\binom{n}{2}-4\end{math} edges.
By Lemma~\ref{lem:decC}, the number of such graphs is the same as the
number of graphs $S$ with 4 edges and without isolated vertices whose
complement \begin{math}\overline{S}\end{math} admits a strictly matched involution.
In the table of Figure~\ref{fig:minusfour} all graphs with 4 edges and
without isolated vertices and their complements are depicted.

It is easy to see that the 6 graph complements 
\begin{math}\overline{S_{1}}\end{math},
\begin{math}\overline{S_{4}}\end{math},
\begin{math}\overline{S_{6}}\end{math},
\begin{math}\overline{S_{7}}\end{math},
\begin{math}\overline{S_{10}}\end{math}, and
\begin{math}\overline{S_{11}}\end{math}
admit a strictly matched involution with the matching edges drawn
vertically in bold lines.

So we only have to prove that the other 5 graph complements
do not admit a strictly matched involution:
\begin{itemize}
\item Consider \begin{math}\overline{S_{2}}\end{math} first. Since none of the edges of
\begin{math}\overline{S_{2}}\end{math} 
has vertices of equal degree and the vertices of the graph
are not a clique, by Lemma~\ref{lem:degree} the graph \begin{math}\overline{S_{2}}\end{math} 
cannot admit
a strictly matched involution.
\item The same argument applies to \begin{math}\overline{S_{3}}\end{math}.
\item In \begin{math}\overline{S_{5}}\end{math}, by Lemma~\ref{lem:degree} the only possible
matching edges would be \begin{math}\{1,4\}\end{math}, \begin{math}\{1,5\}\end{math}, and \begin{math}\{4,5\}\end{math}, and, since
these edges form a triangle, only one of them is possible. But in all
three cases the three vertices not contained in the edge are not a clique.
Therefore no strictly matched involution is possible for~\begin{math}\overline{S_{5}}\end{math}.
\item If we assume that \begin{math}\overline{S_{8}}\end{math} admits
a strictly matched involution~$\sigma$, by Lemma~\ref{lem:degree} the only
possible matching edges are between the vertices of degree 4, namely
between the vertices 2, 3, 4, 5, and 6. Thus there are at most 2 matching
edges. In order to make vertex 1 part of a clique, the clique must consist
of the vertices 1 and 2, or of 1 and 6. By symmetry, we may assume
that the clique consists of the vertices 1 and 2. But there is the edge
\begin{math}\{1,6\}\end{math}, whereas the edge \begin{math}\{1,\sigma(6)\}\end{math} does not exist, which is
a contradiction to the structural theorem (Theorem~\ref{thm:structure}).
\item If we assume that \begin{math}\overline{S_9}\end{math} admits a strictly
matched involution, by Lemma~\ref{lem:degree} the only
possible matching edges are at most two edges on the vertex set \begin{math}\{1,4,5,6\}\end{math}.
However, the vertices 2 and 3 are not connected by an edge, thus
not part of a clique, which is a contradiction.
\end{itemize} 
Note that the graphs \begin{math}S_i\end{math}, \begin{math}1\leq i\leq 11\end{math}, do not have more than
8 vertices.
Thus we have proved that, for every \begin{math}n\geq8\end{math},
\begin{displaymath}C\left(n,\binom{n}{2}-4\right)=6,\end{displaymath}
which proves (\ref{oCc}).
\end{proof}

\begin{remark} From the proof of Theorem~\ref{obs:C}~(\ref{oCc}) we may also 
deduce the following trivial values:
As only \begin{math}\overline{S_{11}}\end{math} has 8 vertices, we conclude that
\begin{displaymath}C\left(7,\binom{7}{2}-4\right)=5.\end{displaymath}
As there are only two of the graphs \begin{math}\overline{S_i}\end{math} 
that admit a strictly matched involution
with more than 6 vertices, we have
\begin{displaymath}C\left(6,\binom{6}{2}-4\right)=4;\end{displaymath}
and by similar arguments we conclude
\begin{displaymath}C\left(5,\binom{5}{2}-4\right)=3\end{displaymath}
and
\begin{displaymath}C\left(4,\binom{4}{2}-4\right)=1.\end{displaymath}
\end{remark}

\section{Open problems and conjectures}\label{sec:open}

The main open problem arising from this paper is the following.

\begin{problem}\label{prob:mainprob}
Determine \begin{math}C(n,m)\end{math} for further values of $n$ and $m$.
\end{problem}

Problem~\ref{prob:mainprob} 
might be attacked by a combinatorial approach, \textit{e.g.}, based on
the decomposition lemma (Lemma~\ref{lem:decC}), or by computational
enumeration approaches.

Even the complexity status of our counting problems are still open, \textit{i.e.},
the following problems might be of further interest.

\begin{problem}
Determine the complexity of Problem~\ref{prob:countA}.
\end{problem}

\begin{problem}
Determine the complexity of Problem~\ref{prob:countB}.
\end{problem}

\begin{problem}
Determine the complexity of Problem~\ref{prob:countC}.
\end{problem}

In Figure~\ref{fig:minusfour} we listed some graphs and their complements.
We observe that some graphs and their complements admit a strictly
matched involution (only $S_1$ and $S_{11}$ and their complements),
some graphs and their complements do not admit a strictly matched involution
(\textit{e.g.}, $S_8$ and its complement), and there are graphs that admit a strictly
matched involution but their complements do not admit a strictly matched
involution (\textit{e.g.}, $S_2$). Thus the following might be interesting.

\begin{problem}
Characterise (or determine the complexity of characterising) 
the class of graphs $G$, for which neither $G$ nor its
complement \begin{math}\overline{G}\end{math} admits a strictly matched involution.
\end{problem}

\begin{problem}\label{prob:last}
Characterise (or determine the complexity of characterising) 
the class of graphs $G$, for which both $G$ and its
complement \begin{math}\overline{G}\end{math} admit a strictly matched involution.
\end{problem}

\appendix

\section{Python3 enumeration program}
\label{sec:append}

\subsection{Idea of the Python3 program}
\label{sec:appendA}

For fixed $n$ and $k$, the enumeration program creates a list of graphs 
admitting a strictly matched involution on $n$ vertices with respect
to some strictly matched involution with
$k$ matching edges according to the structure theorem (Theorem~\ref{thm:structure}) in the following way:
\begin{enumerate}[{\bfseries(Step 1)}]
\item \label{StepA}Based on a list of all graphs $H$ on $k$ vertices, 
for every such graph~$H$, we construct a
graph with two isomorphic copies of $H$, connected by matching edges, and
a clique of size \begin{math}n-2k\end{math}.
\item \label{StepB}For every graph $G$ from the list of graphs from 
Step~\ref{StepA} and for every way to insert
pairs of clique-connecting edges according to part (\ref{oberneuC}) in the proof of 
Theorem~\ref{oberneu}, we construct a graph with the edges of $G$ and
these pairs of clique-connecting edges.
\item \label{StepC}For every graph $G$ from the list of graphs from 
Step~\ref{StepB} and for every way to insert
pairs of cross edges according to part (\ref{oberneuB}) in the proof of 
Theorem~\ref{oberneu}, we construct a graph with the edges of $G$ and
these pairs of cross edges.
\end{enumerate}
Finally, we iterate over all $k$ (\begin{math}k\in\left\{0,1,\ldots,\nhalb\right\}\end{math}) 
and apply isomorphism tests to the union of the 
lists obtained from Step~\ref{StepC} (for all $k$)
to get the list of all (pairwise nonisomorphic) graphs
admitting a strictly matched involution with $n$ vertices (and separate
this list into lists of graphs admitting a strictly matched involution
with $n$ vertices and $m$ edges for \begin{math}m\in\{0,\ldots,\binom{n}{2}\}\end{math};
the number of these graphs is \begin{math}C(n,m)\end{math}).

\subsection{Code of the Python3 program}
\label{sec:appendB}

\begin{verbatim}
import sys
import networkx as nx
import networkx.algorithms.isomorphism as iso

#check whether graph G is isomorphic to a graph of a list
def isIsoTo(G,listOfGraphs):
    for H in listOfGraphs:
        if iso.is_isomorphic(G,H):
            return True
    return False

#sort out isomorphic copies from a list of graphs, 
#return a reduced list
def reduceList(listOfGraphs):
    listOfNonisomorphicGraphs = []
    for G in listOfGraphs:
        isomorphic = False
        for H in listOfNonisomorphicGraphs:
            if iso.is_isomorphic(G,H):
                isomorphic = True
                break
        if not isomorphic:
            listOfNonisomorphicGraphs.append(G)  
    return listOfNonisomorphicGraphs

#creates a graph with n isolated vertices
def edgelessGraph(n):
    G = nx.Graph()
    for i in range(n):
        G.add_node(i)
    return G

#initialises empty list of lists of graphs with n vertices
listOfAllGraphs = {0:[]} 

#creates list of lists of graphs with n vertices
def createListOfAllGraphs(n):
    candidates = []
    candidates.append(edgelessGraph(n))
    for i in range(n-1):
        for j in range(i+1,n):
            newItems = []
            for G in candidates:
                H = G.copy()
                H.add_edge(i,j)
                newItems.append(H)
            for G in newItems:
                candidates.append(G)
    listOfAllGraphs[n] = reduceList(candidates)
    print("Created list of graphs with ",n,
          " vertices. Number: ",len(listOfAllGraphs[n]))

#enumerate all graphs admitting a strictly matched involution
#on n vertices with k matching edges
#do not yet sort out isomorphic copies
def enumerateGraphsWithSMI(n, k):
    #(Step 1) initialise the list 
    #of all graphs admitting a strictly matched involution
    #on n vertices with k matching edges by the empty list
    listOfGraphsWithSMI = []
    #initialise provisional list for these graphs
    newItems = []
    for N in listOfAllGraphs[k]:
        #create union of kK2 and Ks (s=n-2k), which is a subgraph
        #of every graph admitting a strictly matched involution
        #on n vertices with k matching edges
        #indices of vertices of Ks: 0,...,s-1
        #indices of vertices of kK2: s,...,n-1 
        #matching edge (n-2*v-2,n-2*v-1) for v=0,...,k-1
        H=edgelessGraph(n)
        s = n-2*k
        #create Ks --- insert "clique edges"
        for v in range(s):
            for b in range(v+1, s):
                H.add_edge(v, b)  
        #create kK2 --- insert "matching edges"
        for v in range(k):
            H.add_edge(n-2*v-2,n-2*v-1)

        #create edges in the two isomorphic copies of the graph N
        # --- insert pairs of "subgraph edges"
        for v in N.edges():
            H.add_edges_from([(n-2*v[0]-1,n-2*v[1]-1),
                             (n-2*v[0]-2,n-2*v[1]-2)])
        #add the created graph to the provisional list
        newItems.append(H)

    #(Step 2) Now create "clique-connecting edges" by bubbling(n,k)
    #a pair of clique-connecting edges connects 
    #the ends of a matching edge to a clique vertex
    #create all binom(P,m)=binom(2^k+n-2k-1,n-2k) combinations of
    #multisets of subsets of matching edges 
    #on the n-2k clique vertices

    #(Step 2.1): calculate P and m
    twoToThePowerOfK=1
    for i in range(k):
        twoToThePowerOfK=twoToThePowerOfK*2
    m=n-2*k
    P=twoToThePowerOfK+m-1

    #(Step 2.2): initialise 
    #A: incidence array of big set of order P
    #B: multiplicity array
    #start with the first possible values for A
    A=[]
    for i in range(P):
        A.append(0)
    for i in range(m):
        A[i]=1
    B=[]
    for i in range(m):
        B.append(0)

    #(Step 2.3): initialise bubbling stop variable
    bubblingsend=0

    #(Step 2.4): big bubbling loop
    while(bubblingsend==0):

        #(Step 2.4a): create the actual graphs H
        for G in newItems:
            H = G.copy()
            numberZ=0
            cliquecounter=0
            for i in range(P):
                if(A[i]==0):
                    numberZ+=1
                if(A[i]==1):
                    B[cliquecounter]=numberZ
                    cliquecounter+=1
            for cliquecounter in range(m):
                numberZ=B[cliquecounter]
                twoToThePowerOfX=twoToThePowerOfK//2
                for vv in range(k):
                    if(numberZ//twoToThePowerOfX>0):
                        #add pair of "clique-connecting edges"
                        H.add_edges_from([(cliquecounter,n-2*vv-1),
                                          (cliquecounter,n-2*vv-2)])
                        numberZ=numberZ-twoToThePowerOfX
                    twoToThePowerOfX=twoToThePowerOfX//2 
            listOfGraphsWithSMI.append(H)

        #(Step 2.4b): technical bubbling
        leadpointer=-1
        pointercount=0
        Aindex=0
        while((Aindex<P-1) and (leadpointer<0)):
            if(A[Aindex]==1):
                pointercount+=1
                if(A[Aindex+1]==0):
                    leadpointer=Aindex
            Aindex+=1
        if(leadpointer<0):
            bubblingsend=1
        if(bubblingsend==0):
            A[leadpointer]=0
            A[leadpointer+1]=1
            for smaller in range(pointercount-1):
                A[smaller]=1
            for middle in range(pointercount-1,leadpointer):
                A[middle]=0
    #end of bubbling(n,k)
    
    #(Step 3) add "cross edges" between pairs of matching edges
    for i in range(k-1):
        for v in range(i+1,k):
            newItems = []
            for G in listOfGraphsWithSMI:
                H = G.copy()
                H.add_edges_from([(n-2*i-1,n-2*v-2),
                                  (n-2*i-2,n-2*v-1)])
                newItems.append(H)
            #listOfGraphsWithSMI += newItems
            for G in newItems:
                listOfGraphsWithSMI.append(G)
    #return final list
    return listOfGraphsWithSMI

#start main program
#define lower and upper bound xxx of calculation
xxxLower = 0
xxxUpper = 11

#create list of all graphs with <=n/2 vertices
for n in range(xxxUpper//2+1):
    createListOfAllGraphs(n)
print("Number A[n] of graphs admitting a ",
      "strictly matched involution with n vertices:")

#create list of all graphs with SMI and test 
#sort out isomorphic copies
for n in range(xxxLower,xxxUpper+1):
    listOfGraphsWithSMI = []
    for y in range(n*(n-1)//2+1):
        listOfGraphsWithSMI.append([])
    for k in range(n//2+1):
        temporaryList = enumerateGraphsWithSMI(n,k)
        newList = []
        for y in range(n*(n-1)//2+1):
            newList.append([])
        for GG in temporaryList:
            newList[GG.number_of_edges()].append(GG)
        for y in range(n*(n-1)//2+1):
            for GG in newList[y]:
                if not isIsoTo(GG,listOfGraphsWithSMI[y]):
                    listOfGraphsWithSMI[y].append(GG)
    sum=0
    for y in range(n*(n-1)//2+1):
        sum += len(listOfGraphsWithSMI[y])
    print("n =",n," -->  A[n] =",sum)
    print("Results by number of edges:")
    for y in range(n*(n-1)//2+1):
        print(len(listOfGraphsWithSMI[y])," ",end='')
    print()
\end{verbatim}

\nocite{*}
\bibliographystyle{abbrvnat}
\bibliography{extended_gsmi}
\label{sec:biblio}

\end{document}